\documentclass[final,a4paper]{siamltex}
\usepackage{amsmath}
\usepackage{amsfonts}

\newcommand{\Er}{\mathbb{R}}

\newcommand{\calb}{\mathcal{B}}
\newcommand{\calc}{\mathcal{C}}
\newcommand{\cald}{\mathcal{D}}
\newcommand{\calf}{\mathcal{F}}

\newcommand{\calo}{\mathcal{O}}
\newcommand{\calp}{\mathcal{P}}
\newcommand{\calu}{\mathcal{U}}

\newcommand{\calw}{\mathcal{W}}
\newcommand{\ve}{\varepsilon}
\newcommand{\inter}[1]{\operatorname{int}(#1)}
\newcommand{\dist}[1]{\operatorname{dist}#1}
\newcommand{\esssup}[1]{\operatornamewithlimits{ess\,sup}#1}
\newcommand{\norm}[2][{}]{\lVert#2\rVert_{{#1}}}
\newcommand{\vnorm}[2][{}]{\left\Vert#2\right\Vert_{#1}}
\newcommand{\abs}[2][{}]{\lvert#2\rvert_{#1}}
\newcommand{\vabs}[2][{}]{\left\vert#2\right\vert_{#1}}
\newcommand{\one}{\boldsymbol{1}}
\renewcommand{\hat}{\widehat}
\renewcommand{\tilde}{\widetilde}
\newcommand{\tb}{\overline{t}}
\newcommand{\tab}{\overline{\tau}}

\title{A numerical scheme for impact problems}

\author{Laetitia Paoli \thanks UMR 5585 CNRS Analyse Num\'erique,
Facult\'e des Sciences, Universit\'e Jean Monnet, 23 Rue du Docteur
Paul Michelon, 42023 Saint-Etienne Cedex 2, France, 
({\tt paoli@anumsun1.univ-st-etienne.fr})
\and Michelle Schatzman \thanks UMR 5585 CNRS Analyse Num\'erique,
Universit\'e Lyon 1, 69622 Villeurbanne Cedex,
France({\tt schatz@maply.univ-lyon1.fr})
}

\begin{document}
\maketitle

\begin{abstract}
We consider a mechanical system with impact and $n$ degrees of
freedom, written in generalized coordinates. The system is not
necessarily Lagrangian.
The representative point is subject to a constraint: it must stay
inside a closed set $K$ with boundary of class $C^3$.
 We assume that, at
impact, the tangential component of the impulsion is conserved, while
its normal coordinate is reflected and multiplied by a given
coefficient of restitution $e\in[0,1]$: the mechanically relevant
notion of orthogonality is defined in terms of the local metric for
the impulsions (local cotangent metric). We define a numerical scheme
which enables us to approximate the solutions of the Cauchy problem:
this is an \textsl{ad hoc} scheme which does not require a systematic
search for the times of impact. We prove the convergence of this
numerical scheme to a solution, which yields also an
existence result. Without any \textsl{a priori} estimates, the
convergence and the existence are local; with some \textsl{a priori}
estimates, the convergence and the existence are proved on intervals
depending exclusively on these estimates. The technique of proof uses
a localization of the scheme close to the boundary of $K$; this idea
is classical for a differential system studied in the framework of
flows of a
vector field; it is much more difficult to implement here, because
finite differences schemes are only approximately local: straightening
the boundary creates quadratic terms which cause all the difficulties
of the proof.
\end{abstract}



\pagestyle{myheadings}
\thispagestyle{plain}
\markboth{L. PAOLI AND M. SCHATZMAN}{SCHEME FOR IMPACT}

\section{Introduction}

We study in this article a numerical approximation of dynamics with
impact with a finite number of degrees of freedom and a smooth
constraint. 

The set of constraints is denoted $K$ and satisfies the following assumptions:
\begin{subequations}\label{hypK}
\begin{align}
&\text{$K$ is a closed subset of $\Er^d$ with non empty interior;}\\
&\begin{cases}
\text{the boundary $\partial K$ of $K$ is an embedded sub-manifold}\\
\text{of class $C^3$ of $\Er^d$;}
\end{cases}
\\ 
&\text{$K$ lies on only one side of $\partial K$.}
\end{align}
\end{subequations}
It is possible to find a function $\phi$ of class $C^3$ such that
\begin{equation*}
K=\{u\in \Er^d: \phi(u)\ge 0\}
\end{equation*}
and the differential $d\phi$ does not vanish on $\partial
K=\bigl\{u\in \Er^d:\phi(u)=0\bigr\}$. 

Let $f$ be a continuous function
from $[0,T]\times \Er^d\times \Er^d$ to $\Er^d$
which is locally Lipschitz continuous with respect to its last two
arguments, 
and let $M(u)$ be the mass matrix: $u\mapsto M(u)$ is a mapping
of class $C^3$ from $\Er^d$ to the set of symmetric positive definite
matrices. 

The free dynamics of the system are written in generalized coordinates
as
\begin{equation}
M(u)\ddot u = f(\cdot, u, p),\quad  p =M(u)\dot u.\label{eq:1}
\end{equation}

This system is more general than the system obtained in Lagrangian
mechanics, since we want to include possible dissipative terms in the
dynamics of the problem under discussion.

Let us give the few geometric notations which are absolutely necessary
here, since we use a Riemannian metric; the cotangent bundle
$T^*\Er^d$ is identified to $\Er^d \times \Er^d$,
and its elements are denoted as pairs $(u,\xi)$;
at each point $u$ of $\Er^d$ the metric tensor for tangent vectors is
defined by the 
matrix $M(u)$, and the metric tensor for cotangent vectors is defined
by the matrix $M(u)^{-1}$. The scalar product of two vectors $x$ and
$y$ in the tangent space at
$u$ is denoted by $\langle x, y\rangle_u$; coordinate-wise it can be
expressed as $x^T M(u)y$ where $x$ and $y$ are column vectors. The
scalar product of 
two vectors $\xi$ and $\eta$ in the
cotangent space at $u$ is denoted by $\langle\xi,\eta \rangle_u^*$ and
coordinate-wise it is equal to $\xi^T M(u)^{-1}\eta$. The corresponding
norms of vectors and covectors are denoted respectively by $\abs[u]{x}$
and $\abs[u]{\xi}^*$.

Therefore, a cotangent vector $(u,\xi)$
belonging to $T^*\Er^d$ is orthogonal to the cotangent vector
$(u,\eta)$ iff $\langle \xi, \eta\rangle_u^*$ vanishes. 

With these notations,
if the velocity of the system is $\dot u$, the generalized impulsion
is $M(u)\dot u=p$ and $(u,p)$ belongs to the cotangent space
$T^*\Er^d$. Whenever we take the orthogonal of a vector or
a vector subspace of the tangent or the cotangent space at $u$, we
always use the relevant metric tensor; therefore it is important to
know which of the vectors under consideration are cotangent and which
are tangent. Of course, all the differential forms are cotangent vectors.

Let us describe now the system satisfied by the problem with impact:
we replace~\eqref{eq:1} by
\begin{equation}
M(u)\ddot u =\mu + f(\cdot, u, p),\label{eq:2}
\end{equation}
and since we cannot expect to have global solutions in general,
$\mu$ is an unknown measure on $[t_0,t_0+\overline{\tau}]$ with values in
$\Er^d$ which
describes the reaction of the constraints: $\mu$ has 
the following properties: if $d\phi$ denotes the differential of
$\phi$, then
\begin{subequations}
\begin{align}
&\supp(\mu)\subset \{t\in [t_0,t_0+\overline{\tau}]: \phi(u(t))=0\},\label{eq:3}\\
& \mu = \lambda d\phi(u),\label{eq:4}\\
&\lambda \ge 0 \text{ almost everywhere on $[t_0,t_0+\overline{\tau}]$}.\label{eq:5}
\end{align}
\end{subequations}

We require the following functional properties for $u$:
\begin{subequations}
\begin{align}
&\begin{cases}
\text{$u$ is a continuous function taking its values in $K$}\\
\text{for all $t\in [t_0,t_0+\overline{\tau}]$,}
\end{cases}\label{eq:7}
\\
&\text{$\dot u$ is of bounded variation over
$[t_0,t_0+\overline{\tau}]$.}\label{eq:8}  
\end{align}
\end{subequations}

If $\dot u$ is of bounded variation, $p$ is also of bounded variation.
Assume that $u(t)$ belongs to $\partial K$; we
decompose $p(t-0)$ and $p(t+0)$ on $\Er d\phi(u(t))\oplus
d \phi(u(t))^\perp$; here the $\perp$ sign means the orthogonality
with respect to the local cotangent metric.
We integrate~\eqref{eq:2} on a small
neighborhood of $t$, relation~\eqref{eq:4} implies that the component
of $p(t-0)$ on $d \phi(u(t))^\perp$ is conserved.

Therefore, we have to make a supplementary assumption in order to have
a complete description of the impact; we choose a constitutive law of
the impact using a coefficient of restitution: thus we will assume
that there exists $e\in [0,1]$ such that the component of $\dot
p(t+0)$ along $\Er d\phi(u)$ is equal to $-e$ times the component of
$p(t-0)$ on $\Er d\phi(u)$. In other words,
we have
\begin{equation}
p(t+0)= p(t-0)-(1+e)\frac{\langle
d\phi(u(t)),p(t-0)\rangle_{u(t)}^*}{\langle
d\phi(u(t)),d\phi(u(t))\rangle_{u(t)}^*}d\phi(u(t)).
\label{eq:6}
\end{equation}

The set of admissible initial data $\mathbb{D}$ will be
\begin{equation}
\begin{split}
\mathbb{D}&=\Bigl\{(t_0,u_0, p_0)\in [0,T)\times K\times \Er^d:\\
&\text{if $u_0\in 
\partial K$, then 
$\langle p_0, d\phi(u_0)\rangle_{u_0}^*\ge 0$}\Bigr\}.
\end{split}
\label{condinit}
\end{equation}
This choice is equivalent to the convention that there is no impact at
the initial time $t_0$.

Given initial conditions $(t_0,u_0, p_0)\in\mathbb{D}$, we
require that the
following Cauchy data be satisfied:
\begin{equation}
u(t_0)=u_0,\label{position}
\end{equation}
and
\begin{equation}
p(t_0)= p_0.\label{eq:9}
\end{equation}

For all initial data $(t_0,u_0, p_0)\in \mathbb{D}$ we  will obtain the
existence of a local solution
to~\eqref{eq:2},~\eqref{eq:3},~\eqref{eq:4},~\eqref{eq:5}
and~\eqref{eq:6}  belonging
to the functional class defined by~\eqref{eq:7} and~\eqref{eq:8} and
satisfying the initial conditions~\eqref{position} and~\eqref{eq:9}.

The existence of this local solution is obtained by defining a
numerical scheme, whose convergence will
be shown in appropriate functional spaces; the limit of the
approximation will be 
a solution of our problem.

The distance on $\Er^d$ is defined with the help of the
Riemannian metric: if $s\mapsto u(s)$ is a $C^1$ mapping from $[a,b]$
to $\Er^d$, the Riemannian length of the image of $u$ is
\begin{equation*}
\ell(u)=\int_a^b \vabs[u(s)]{\dot u(s)}\, ds.
\end{equation*}
This curve length is invariant by a diffeomorphic change of
parameter.
Therefore, we may assume that $a=0$ and $b=1$.
The distance from $x$ to $y$ is the lower bound of the length of the
curves from $x$ to $y$, or in other words:
\begin{equation*}
\dist(x,y)= \inf\{\ell(u):\quad u\in C^1([0,1]), \quad u(0)=x,\quad
u(1)=y\}. 
\end{equation*}
It is classical that the lower bound is attained on the
geodesics for the given Riemannian metric; it is also known that for
each point $x$ there exists $r>0$ such that if $\dist(x,y)\le r$ there
is only one geodesic from $x$ to $y$.

We denote  by $\dist(x,
E)$ the Riemannian distance of a point $x$ to a set $E$.

Under assumptions \eqref{hypK}, a projection on $\partial K$ can be
defined uniquely on an appropriate neighborhood of $\partial K$; more
precisely, 
for all compact $\calc\subset \partial K$, there exists a neighborhood 
of $\calc$ on which the projection $P_{\partial K}$ is uniquely
defined, and there exists a unique geodesic joining 
a point of this neighborhood to its
projection. This projection $P_{\partial K}$ is 
characterized by the relation
\begin{equation}
\forall y\in \partial K, \quad \dist(P_{\partial K} x, x)\le\dist(y,x).
\end{equation}
This projection is of class $C^2$.

For all $x$ in $\partial K$, denote by $N(x)$ the interior unit normal
vector: this means that $\abs[x]{N(x)}$ is equal to $1$ and that it
is orthogonal to the 
tangent space at $P_{\partial K}x$ with respect to the scalar product
in the tangent space, i.e. for all $y$ such that $d\phi(x) y$
vanishes, $\langle y, N(x)\rangle_x=0$. The smoothness of $\partial K$
implies that the mapping $z\mapsto 
N(z)$ is 
of class $C^2$.

When the geodesic from $x$ to
$P_{\partial K} x$ is unique it is tangent at $P_{\partial K} x$ to
$N(P_{\partial K}x)$.

Starting from this projection on $\partial K$, we can define a projection on
$K$ as follows: for each compact $\calc$ included in $K$, there exists
a relatively compact neighborhood $\calu$ of $\calc$ on which $P_K$ is defined by
\begin{equation}
P_K(x)=
\begin{cases}
P_{\partial K}(x) & \text{if $ x\notin K$}\\
x&\text{otherwise.}
\end{cases}
\end{equation}

The reader will check that $P_K$ is Lipschitz continuous over
$\calu$ and that $P_Kx$ realizes the minimum of
the distance from $x$ to $K$.

Given two positive numbers $h^*\le 1$ and $T$, assume that
$F$ is a continuous function from $[0,T]\times \Er^d\times \Er^d
\times \Er^d 
\times [0, h^*]$ to $\Er^d$, which is locally Lipschitz continuous with
respect to its 
second, third and fourth arguments; assume moreover that $F$ is
consistent with $f$, i.e. that for all $t\in [0,T]$,
for all $u$ and $v$ in $\Er^d$
\begin{equation}
F(t,u,u,v,0)=M(u)^{-1}f(t,u,M(u)v).\label{consistance}
\end{equation}

We approximate the solution of
\eqref{eq:2}, \eqref{eq:3}, \eqref{eq:4}, \eqref{eq:5},
\eqref{eq:7}, \eqref{eq:8}, \eqref{position},
\eqref{eq:9} by the following numerical scheme:
the initial values $U^0 $ and $U^1$ are given by the initial position
\begin{align}
U^0&=u_0,\label{u0}\\
\intertext{and the position at the first time step}
U^1&=u_0 + hM(u_0)^{-1}p_0 +hz(h),\label{u1} 
\end{align}
where $z(h)$ tends to $0$ as $h$ tends to $0$.

We will use systematically henceforth the notation
\begin{equation}
t_m =t_0+mh.\label{eq:74}
\end{equation}

Given $U^{m-1}$ and $U^m$, $U^{m+1}$ is defined by the relations
\begin{align}
U^{m+1}&=-eU^{m-1}+(1+e)P_K\left(\frac{2U^m -(1-e)U^{m-1}+ h^2
F^m}{1+e}\right)\label{un}\\
\intertext{and}
F^m&=F\left(t_m,U^m,U^{m-1},\frac{U^{m+1}-U^{m-1}}{2h},h\right)\label{fn}
\end{align}
provided that $U^{m+1}$ is unique in a neighborhood of $U^m$.

A commentary on the construction of this scheme from the point of view
of convex 
analysis will be useful here. We refer to the book of Rockafellar
\cite{rock} for 
more information on the basic ideas in convex analysis to be used below.

Let us assume provisionally that the set of constraints
$K$ is convex and that the mass matrix is equal to the
identity matrix on $\Er^d$. Then the Riemannian structure of $\Er^d$
is simply its Euclidean structure. 

Recall that the indicator function $\psi_K$ of a closed convex set $K$ is defined
by
\begin{equation}
\psi_K(x)=
\begin{cases}
0 &\text{if $x \in K$},\\
+\infty &\text{otherwise},
\end{cases}
\end{equation}
and its sub-differential $\partial \psi_K$ is a function from $K$ to the set of
closed convex sets given by
\begin{equation}
\partial \psi_K(x)=
\begin{cases}
\{0\} &\text{if $x\in \inter{K}$},\\
\Er^- N(x) & \text{if $x\in \partial K$}.
\end{cases}
\end{equation}
For all $\lambda >0$, the multivalued equation
\begin{equation}
x+ \lambda \partial \psi_K(x)\ni f \label{multi1}
\end{equation}
has a unique solution given by
\begin{equation}
x = P_K(f), \label{multi2}
\end{equation}
where $P_K$ is the usual projection on the closed convex set $K$ in
Euclidean $\Er^d$.

In the annoucement \cite{crasps}, where we assumed that the set of
constraints $K$ was convex and the geometry was Euclidean, we had
defined the numerical scheme 
by the multivalued equation
\begin{equation}
\frac{U^{m+1} - 2 U^m + U^{m-1}}{h^2} + \partial \psi_K\left(\frac{U^{m+1} + e
U^{m-1}}{1+e}\right)\ni F^m. \label{original}
\end{equation}
We may rewrite \eqref{original} as
\begin{equation}
\begin{split}
&\frac{U^{m+1} + e U^{m-1}}{1+e} + \frac{h^2}{1+e} \partial \psi_K
\left(\frac{U^{m+1} + e U^{m-1}}{1+e}\right) \\
&\quad \ni \frac{2U^m - (1-e)U^{m-1} + h^2
F^m}{1+e},
\end{split}
\end{equation}
which reduces, thanks to \eqref{multi1} and \eqref{multi2} to relation
\eqref{un}. 

If we generalize \eqref{original} to a non convex $K$ with a general mass
matrix, we cannot use the apparatus 
of convex analysis, and there is no good reason to use the even more
technical 
apparatus of non-convex analysis \textsl{\`a la } Clarke: this theory
is useful when the corners of $K$ are not convex; in the mechanical
setting, corners are convex, since they appear
as the intersection of smooth sets of constraints. Here, the problem
is even simpler because we do not have any corners.

The boundary $\partial K$ is
smooth, and as we expect that for small $h$, the $U^m$'s will stay
close to $K$, 
we still have a projection of $(2U^m - (1-e)U^{m-1} + h^2
F^m)/(1+e)$ on $K$, and thus we start from \eqref{un} to define the
numerical scheme. 

The original definition reappears as follows: define
\begin{equation}
W^m= \frac{2U^m - (1-e) U^{m-1} + h^2 F^m}{1+e},\label{Wn}
\end{equation}
that will be used in many places in the upcoming proofs.
With this definition, \eqref{un} is rewritten as
\begin{equation*}
U^{m+1} = - eU^{m-1} + (1+e) P_K(W^m).
\end{equation*}
Hence, if we define
\begin{equation}
Z^m = \frac{U^{m+1} + e U^{m-1}}{1+e},\label{Zn}
\end{equation}
we find that 
\begin{equation}
Z^m= P_K(W^m).\label{eq:78}
\end{equation}
If we subtract \eqref{Zn} from \eqref{Wn}, we can see that
\begin{equation}
\frac{U^{m+1} - 2 U^m + U^{m-1}}{h^2} + \frac{1+e}{h^2} (W^m - Z^m)= F^m
\label{poriginal}
\end{equation}
which reduces to \eqref{original} in the convex case with a trivial
mass matrix.

Another way of writing \eqref{poriginal} is to define the discrete velocity $V^m$
by
\begin{equation}
V^m= \frac{U^{m+1} - U^m}{h} \label{Vn}.
\end{equation}
Then, \eqref{poriginal} can be rewritten as
\begin{equation}
V^m - V^{m-1} - h F^m = \frac{(1+e)(Z^m - W^m)}{h}.\label{diffVn}
\end{equation}

A strict contraction argument in $\Er^d$ gives the existence of a unique
$U^m$ for small values of $m$ and $h$. As the projection on $K$ is uniquely
defined only in a neighborhood of $K$, and is only Lipschitz continuous, the
iteration of a fixed point argument might request smaller and smaller bounds
on the time step $h$, and there is no guarantee that we could
integrate numerically 
on a time interval bounded from below, for any initial time step size.

It should be noted that this difficulty is specific to the non convex
case.

Let us outline now the structure of the article and of the proofs. In
the one-dimensional case, the main estimates are given by
lemma~\ref{monodim}, in section~\ref{sec:heart-estimates}. 
In section~\ref{sec:Existence}, we will straighten the boundary, a
natural geometrical idea. 

While the
system~\eqref{eq:2}--\eqref{eq:6} is nicely transformed under a
diffeomorphism, the numerical scheme~\eqref{u0}, \eqref{u1},
\eqref{un} and~\eqref{fn} does not behave well under
diffeomorphism. The reason is that a numerical scheme is not a local
object: when we define a discrete velocity by subtracting $U^m$ from
$U^{m+1}$, we use locally a vector structure which is not 
intrinsic from the point of
view of differential geometry. In particular, if we apply a diffeomorphism to the
numerical scheme, we will find another numerical scheme which will
look much more complicated than the previous one, since it will
contain a number of small term which show the lack of an intrinsic
description of the scheme. After a very technical proof, we find two
constants $C_3$ and $\tau$ such that for initial data in a compact
subset of the admissible set, and for all small enough $h$ and all
$m\le \tau/h$, the discrete velocity is bounded:
\begin{equation*}
\sup\vabs{V^m}\le C_3.
\end{equation*}

Since uniqueness is not true in general \cite{bressan1},
\cite{schatz78}, and hypotheses of analyticity are often but not
always used for the
proof of uniqueness \cite{percivaleuniq}, \cite{percivaleuniq2},
\cite{schatz98}, \cite{carpasc}, \cite{ballard}, the proof of 
convergence of the numerical approximation is delicate also for this reason.

However, there is a bonus: all the effort made to prove the local
convergence of the numerical scheme provides us with a local existence
proof for our problem. In
sections~\ref{sec:Estim-accel}, \ref{sec:Variational}, ~\ref{sec:Transmission-energy}
and~\ref{sec:Initial-conditions}, we prove estimates on the discrete
acceleration, we establish the variational properties of the limit of the numerical
scheme, and we study the transmission of energy at impact, as well as
the passage to the limit for the initial conditions. All these results
are obtained under the assumption that on a certain time interval
starting at $t_0$, the discrete velocity is bounded independently of
the time step. 

As a preliminary to the global existence proof, we give \textsl{a priori}
estimates on problem~\eqref{eq:2}--\eqref{eq:9} in
section~\ref{sec:apriori}, which is completely independent from the
remainder of the article.

In section~\ref{sec:Global}, we establish a very weak
semi-continuity for the supremum of the local norm of the discrete
velocities; this result enables us to obtain a global existence and
convergence theorem.

This article is of a theoretical nature: the existence result obtained
here is a generalization of \cite{schatz78}, \cite{buttazzopercivale},
\cite{percivaleweak}, 
 \cite{paolithese}, \cite{pasc93}.

The numerical scheme analyzed here has been implemented in the case of
a trivial mass matrix in \cite{paolithese}, \cite{pasc95}, \cite{papasc}, 
\cite{resopasc}. In all these articles, we compared the
performances of this scheme with those of a method based on the
detection of impact. When the impact times are isolated, the scheme by
detection of impacts is more precise than the present scheme. 
As soon as the restitution coefficient is strictly less than one, we
find systematically non-isolated impact times.
In all cases, the present scheme is substantially faster. Since the
phenomena that we want to approximate are highly nonlinear and
often very sensitive to the initial data, the issue of precision is
not necessarily crucial. Our numerical experiments show that the
performance of the present numerical scheme is quite satisfactory from
the point of view of qualitative conclusions.

The case of a non-trivial mass matrix, and a stiff system, indeed the
case of the discretization of a beam has been adressed in \cite{PS1999}.

Let us remark that many articles have been devoted to the
problem treated here, under the assumption of anelasticity,
i.e. a situation where the normal component of the impulsion vanishes
after the impact; Moreau applied
Gauss' principle of least constraint to unilateral problems in order
to justify his choice of anelastic impact \cite{MR30:5527}, which
eventually led him to sweeping processes \cite{MR58:22889}, followed
by \cite{MR58:30612}, \cite{MR58:30613}; dry friction enters in
Moreau's work as \cite{MR81e:73091}; frictionless anelastic impact starts as
\cite{moreau1}, and the mathematical theory is tackled by
M. Monteiro-Marques in a series of articles: his main contributions
are \cite{MR94g:34003} for the general theory of differential
inclusions, \cite{MR95d:34028} for one-dimensional dynamics with
friction, \cite{MR96k:34007} which adds percussion to the previous
framework; this work is improved as \cite{MR99a:34028}, where dynamics
of $n$ particles on a plane with normal friction are considered. The
discretization approach has been taken up by Monteiro-Marques and
Kuntze in \cite{MR99c:34027}, but most significantly by Stewart and
Trinkle: they use that approach in \cite{MR97b:34013}, \cite{MR97d:70004} and
\cite{MR98i:70009}; the real coronation is the
beautiful and difficult article of Stewart \cite{MR1664526}, which
concludes the study 
of dynamics with friction and anelastic impact for a finite number of
degrees of freedom, and one constraint, and still important results in
the multiple constraint case.

The philosophy of this long list of works is
somewhat different from ours: we feel that not all impacts are anelastic, and
we were originally motivated by continuous media; thus, we wanted
to develop methods which work well for stiff systems of ordinary
differential equations. From this point of view, any method which has
to calculate with some precision the impact times is doomed to
failure. On the other hand, the precision of the method 
presented here needs improvement, and globally, it would make sense to
agree on benchmarks which would enable the end-user to decide between
different numerical methods.

\section{The heart of the estimates}\label{sec:heart-estimates}

In the one-dimensional case, the main estimate on the numerical scheme is
described in the following lemma; we recall the definition
\begin{equation*}
r^+=\max(r,0).
\end{equation*}

\begin{lemma}\label{monodim}Let the real-valued sequence
$\bigl(y^m\bigr)_m$ satisfy the following 
recurrence relation for all $m\ge 1$:
\begin{equation}
y^{m+1}=-e y^{m-1} + \bigl(2 y^m -(1-e)y^{m-1}\bigr)^+ +
h^2\lambda^m.\label{eq:106} 
\end{equation}
Then, for all $m\ge 2$, the discrete velocity
\begin{equation}
\eta^m=\bigl(y^{m+1}-y^m\bigr)/h\label{eq:107}
\end{equation}
satisfies the estimate
\begin{equation}
\vabs{\eta^m} \le \max\bigl(\vabs{\eta^{m-1}},e\vabs{\eta^{m-2}}\bigr)+
h\vabs{\lambda^m}+ h\vabs{\lambda^{m-1}}.\label{eq:131}
\end{equation}
\end{lemma}

\begin{proof}
Assume first that $2y^m-(1-e)y^{m-1}$ is non negative, and substitute
$y^{m+1}=y^m + h\eta^m, y^{m-1}= y^m - h\eta^{m-1}$
into~\eqref{eq:106}; we obtain
\begin{equation*}
\eta^{m}=\eta^{m-1} + h \lambda^m,
\end{equation*}
so that
\begin{equation}
\vabs{\eta^{m}}\le \vabs{\eta^{m-1}} + h \vabs{\lambda^m}.\label{eq:102}
\end{equation}
Assume now that $2y^m-(1-e)y^{m-1}$ is strictly negative. On one hand,
\eqref{eq:106} implies the relation
\begin{equation*}
\eta^m = e\eta^{m-1} -\frac{1+e}{h} y^m + h\lambda^m;
\end{equation*}
the assumption on the sign of $2y^m-(1-e)y^{m-1}$ is equivalent to
\begin{equation*}
\frac{(1+e)y^m}{h} < -(1-e)\eta^{m-1},
\end{equation*}
and therefore
\begin{equation}
\eta^m > \eta^{m-1} + h\lambda^{m}.\label{eq:103}
\end{equation}
On the other hand, we subtract from the relation
\begin{align*}
y^{m+1} + ey^{m-1}&= h^2\lambda^m\\
\intertext{the inequality implied by~\eqref{eq:106} with $m$
substituted by $m-1$:}
y^{m}+ ey^{m-2} &\ge h^2\lambda^{m-1},
\end{align*}
and we infer that
\begin{equation}
\eta^m \le - e\eta^{m-2} + h \bigl(\lambda^m -
\lambda^{m-1}\bigr).\label{eq:104}
\end{equation}
When we summarize~\eqref{eq:102},~\eqref{eq:103} and~\eqref{eq:104},
we find~\eqref{eq:131}.
\end{proof}

Later on, we will give a $d$-dimensional version of~\eqref{eq:131},
where the main difference is due to geometric effects: there will be a
term resembling $\lambda^m$, and the game will be to prove a bound on
this term.

\section{Existence of $(U^m)_{0 \le m \le \lfloor \tau/h
\rfloor}$ for some $\tau>0$}\label{sec:Existence}

We use systematically the floor and ceiling notations: when $r$
is a real number, the floor $\lfloor r\rfloor$ of $r$ is the largest
integer at most equal to $r$, and the ceiling $\lceil r \rceil$ is the
smallest integer at least equal to $r$.

The main result of this section is the existence of a number $\tau>0$
such that for all small enough $h$ and all
$m\le
\lfloor
\tau/h\rfloor$ there exists indeed a discrete solution of \eqref{un} and
\eqref{fn}, whose discrete velocity is bounded independently of
$h$. In fact, we prove a stronger result: provided that the first two
discrete velocities are bounded, we find a uniform lower bound on
$\tau$ when the initial position belongs to a compact subset of $K$.

We prove first the existence of $U^2$ under appropriate
assumptions on $U^0$ and $U^1$. This proof decomposes in two lemmas:
the first lemma is strictly an initial condition statement, in which
no uniformity with respect to initial conditions can be obtained. The
second one will be used in the foregoing induction proofs.

\begin{lemma}\label{implicit0} For all $(t_0, u_0,M(u_0)v_0)\in \mathbb{D}$,
for all $U^1$ satisfying~\eqref{u1}, and for all small enough $h$,
there exists a solution $U^2$ of~\eqref{un} for $m=2$ satisfying
\begin{equation*}
\vabs[u_0]{U^2 -U^1}\le 2\abs[u_0]{v_0}h.
\end{equation*}
\end{lemma}

\begin{proof}
Let $r>0$ be such that $P_K$ is Lipschitz continuous on 
\begin{equation*}
B_{u_0}(u_0,r)=\bigl\{u\in\Er^d: \vabs[u_0]{u-u_0}\le r\bigr\}.
\end{equation*}
Define $\tilde C_1$ by
\begin{equation*}
\begin{split}
\tilde C_1&=\max\bigl\{ \vabs[u_0]{F(t,u,u',0,h)}:\> t\in[0,T],\>
\vabs[u_0]{u-u_0}\le r,\\& \vabs[u_0]{u'-u_0}\le r,\> h\in
[0,h^*]\bigr\},
\end{split}
\end{equation*}
and let $\tilde L$ be the Lipschitz constant defined by
\begin{equation*}
\begin{split}
\tilde
L&=\sup\biggl\{
\frac
{\vabs[u_0]{F(t,u,u',v,h)-F(t,u,u',v',h)}}
{\vabs[u_0]{v-v'}}:
\vabs[u_0]{u-u_0}\le r,\\& \vabs[u_0]{u'-u_0}\le r,\>
\vabs[u_0]{v}\le 2\vabs[u_0]{v_0}+1,\> \vabs[u_0]{v'}\le
2\vabs[u_0]{v_0}+1,\\& v\neq v',\> h\in[0,h^*] 
\biggr\}.
\end{split}
\end{equation*}
Finally, let $\tilde \gamma$ be the Lipschitz constant of $P_K$
defined by
\begin{equation*}
\tilde \gamma=\sup\left\{\frac{\vabs[u_0]{P_K u -
P_Ku'}}{\vabs[u_0]{u-u'}}: \vabs[u_0]{u-u_0}\le r,
\>\vabs[u_0]{u-u_0}\le r,\> u\neq u'
\right\}.
\end{equation*}
There exists a function $\hat z(t)$ which is bounded in a neighborhood
of $0$
such that for small positive values of $t$:
\begin{equation}
P_K(u_0 +tv_0)=u_0 + t v_0 + t^2\hat z(t);\label{eq:16}
\end{equation}
indeed, if $v_0$ vanishes, or if $u_0$ belongs to $\inter K$, or if
$u_0$ belongs to 
$\partial K$ and the scalar product $\langle v_0, N(u_0)\rangle_{u_0}$ is strictly
positive, $\hat z$ vanishes; if $u_0$ belongs to $\partial K$ and
$\langle v_0, N(u_0)\rangle_{u_0}$ vanishes, while $v_0$ does not
vanish,~\eqref{eq:16} is a consequence of the smoothness of
$P_{\partial K}$ in a neighborhood of $u_0$ : for the values of $t$
for which $u_0 +tv_0$ belongs 
to $K$, $\hat z$ vanishes; for the values of $t$ for which $u_0 +
tv_0$ does not belong to $K$, a Taylor expansion shows that
\begin{equation*}
P_{\partial K}(u_0 +tv_0)= u_0 + t v_0 + O(t^2),
\end{equation*}
hence~\eqref{eq:16}.
With the change of variable $U^2=U^1 + tV^1$, equation~\eqref{un} is
equivalent to 
\begin{equation*}
v=\tilde G(v)
\end{equation*}
where the function $\tilde G$ is defined by
\begin{equation*}
\begin{split}
\tilde G(v)&=-V^0 +\frac{1+e}{h}\biggl[P_K\biggl(U^0 + \frac{2h}{1+e} V^0\\
&+\frac{h^2}{1+e}
F\left(t_1,U^1,U^0,\frac{V^0+v}{2},h\right)\biggr)-U^0\biggr].
\end{split}
\end{equation*}
Let us check that $\tilde G$ is a strict contraction on
$B_{u_0}\bigl(0,2\vabs[u_0]{v_0}+1\bigr)$: if  
$\vabs[u_0]{v}\le 2\vabs[u_0]{v_0}$, then
\begin{equation*}
\vabs{\frac{v+V^0}{2}}\le \vabs[u_0]{v_0} +\frac12 \vabs[u_0]{v_0
+z(h)}+\frac12; 
\end{equation*}
therefore, for $h$ small enough, $\vabs[u_0]{(v+V_0)/2}$ is at most
equal to $2\vabs[u_0]{v_0}$, and we can use the definitions of $\tilde
L$ and $\tilde C_1$:
\begin{equation}
\vabs[u_0]{\left(F(t_1,U^1,U^0,\frac{V^0+v}{2},h\right)
}\le \tilde C_1 + \tilde L\bigl(2
\vabs[u_0]{v_0}+1\bigr).\label{eq:48}
\end{equation}
We estimate $G(v)$ as follows: by the triangle inequality, and the
Lipschitz condition on $P_K$,
\begin{equation*}
\begin{split}
&\vabs[u_0]{G(v)} 
\le \frac{1+e}{h}
 \tilde \gamma \vabs[u_0]
{U^0 +
\frac{2h\bigl(V^0-v_0\bigr)}{1+e} + \frac{h^2}{1+e} F - u_0}\\
&\quad+\vabs[u_0]{-V^0
+\frac{1+e}{h}\left(P_K\left[u_0+\frac{2hv_0}{1+e}\right]
-u_0\right)}.
\end{split}
\end{equation*}
We apply~\eqref{eq:16}, \eqref{u1} and~\eqref{eq:48}, and we find
\begin{equation*}
\begin{split}
\vabs[u_0]{G(v)} &\le \tilde \gamma \Bigl[2\vabs[u_0]{z(h)}  +
h\bigl(\tilde C_1 + \tilde L(2\vabs[u_0]{v_0}+1)\bigr)
\Bigr]\\
&\quad +\vabs[u_0]{v_0} +\vabs[u_0]{z(h)}
+\frac{4h}{1+e}\vabs[u_0]{\hat z\left(\frac{2h}{1+e}\right)}.
\end{split}
\end{equation*}
Therefore, for $h$ small enough, $G$ maps
$B_{u_0}\bigl(0,2\vabs[u_0]{v_0}+1\bigr)$ to itself. Morover, the
Lipschitz constant of $G$ on this ball is at most equal to $\tilde
\gamma\tilde L h/2$. This proves that $G$ has a fixed point in
$B_{u_0}\bigl(0,2\vabs[u_0]{v_0}+1\bigr)$ for small enough values of $h$
and completes the proof of the lemma.
\end{proof}

Here is the statement of the uniformizable estimate, which will be
used throughout the induction. We will say that $U^0$ and
$U^1$ satisfy condition $E\bigl(\overline{u},r_0,C_2,h\bigr)$ if
\begin{equation}
\vabs{U^0-\overline{u}}\le r_0, \quad \vabs{U^1-\overline{u}}\le r_0,\quad
\vabs{U^1-U^0}\le C_2 h\label{eq:17}
\end{equation}
and moreover $U^2$ is uniquely defined in $B(\overline{u},r_0)$
by
\begin{equation}
\frac{U^2 + eU^0}{1+e}= P_K\left(\frac{2U^1 -(1-e)U^0 +h^2
F^1}{1+e}\right),\label{eq:21} 
\end{equation}
and the following inequalities are satisfied:
\begin{equation}
\vabs{U^2-\overline{u}}\le r_0, \quad \vabs{U^2-U^1}\le C_2 h.\label{eq:20}
\end{equation}

\begin{lemma}\label{implicit} For all $\overline{u}\in K$, there
exists $r_0$ such that for
all $C_2>0$ it is possible to find $h_1>0$ and $C'_2<\infty$ 
with the following properties:
for all $h\le h_1$ and for all choice of
$U^0$, $U^1$ satisfying condition
$E\bigl(\overline{u},r_0,C_2,h\bigr)$
i.e.~\eqref{eq:17}, \eqref{eq:21} and \eqref{eq:20},
there exists a unique $U^3$
satisfying~\eqref{un} for $m=2$ and the estimate
\begin{equation*}
\vabs{U^3 -U^2}\le C'_2 h.
\end{equation*}
\end{lemma}

\begin{proof}Subtract~\eqref{un} for $m=1$ from~\eqref{un} for $m=2$;
with the change of variable $U^3 =U^2 + hV^2$, we have to show the
existence of a solution of
\begin{equation*}
\begin{split}
V^2 &=-e V^0 +\frac{1+e}{h} P_K\left(\frac{2U^1 + 2h V^1 -(1-e)(U^0
+hV^0) +h^2 F^2}{1+e}\right)\\
&-\frac{1+e}{h}P_K \left(\frac{2U^1  -(1-e)U^0
+h^2 F^1}{1+e}\right).
\end{split}
\end{equation*}
If we denote by $G(V^2)$ the right hand side of the above equation, we
have to choose a parameter $C'_2$ such that $G$ will be a strict
contraction of the ball $B(0,C'_2)$ into itself. Let $r_0$ be such
that $P_K$ is Lipschitz continuous on $B(\overline{u},2r_0)$; denote
by $\gamma$ the Lipschitz constant of $P_K$ over this ball, and define
\begin{equation}
C'_2=(3\gamma +1) C_2;\label{eq:22}
\end{equation}
let $C_1$ be given by
\begin{equation}
\begin{split}
C_1=&\sup\bigl\{\vabs{F(t,u,u',0,h)}: t\in[0,T],\>
\vabs{u-\overline{u}}\le 2 r_0,\\ 
&\qquad\vabs{u'-\overline{u}} \le 2 r_0, \> h\in[0,h^*]\bigr\}.
\end{split}
\label{eq:141}
\end{equation}
Denote finally by $L$ the Lipschitz constant of $F$ defined as
follows:
\begin{equation}
\begin{split}
L &=\sup\biggl\{
\frac{\vabs{F(t,u,u',v,h)-F(t,u,u',v',h)}}{\vabs{v-v'}}: t\in [0,T],\>
\vabs{u-\overline{u}}\le 2r_0,\\&\vabs{u'-\overline{u}}\le 2r_0,
 \vabs{v}\le C'_2,\>
\vabs{v'}\le C'_2,\> h\in[0,h^*],\> v\neq v'\biggr\}.
\end{split}\label{eq:23}
\end{equation}
Then, we have the estimate for $\vabs{v}\le C'_2$:
\begin{equation*}
\begin{split}
\vabs{F\bigl(t_2,U^1,U^2,(v+V^1)/2,h\bigr)}&\le C_1+LC'_2,\\
\vabs{F\bigl(t_1,U^0,U^1,(V^1+V^0)/2,h\bigr)} &\le C_1 + L C'_2.
\end{split}
\end{equation*}
It is straightforward that
\begin{equation*}
\begin{split}
&\max\biggl(\vabs{\frac{2U^1 -(1-e)U^0 +h^2F^1}{1+e}-\overline{u}},\\
&\qquad\vabs{\frac{2U^2 -(1-e)U^1 +h^2F^2}{1+e}-\overline{u}}\biggr)\\
&\quad \le r_0
+\frac{2hC'_2}{1+e}+\frac{h^2(C_1 +LC'_2)}{1+e}.
\end{split}
\end{equation*}
Therefore, if $h_1$ satisfies the estimate
\begin{equation}
\frac{2h_1C'_2}{1+e}+\frac{h_1^2(C_1 +LC'_2)}{1+e} \le r_0,\label{eq:29}
\end{equation}
we may use the Lipschitz continuity of $P_K$ on the ball of radius
$2r_0$ about $\overline{u}$, and we find that if $v$ belongs to $B(0,C'_2)$,
\begin{equation*}
\vabs{G(v)}\le e C_2 +\gamma\bigl((3-e)C_2 + 2h(C_1+LC'_2)\bigr);
\end{equation*}
We observe that $\gamma$ is at least equal to $1$, since $DP_{\partial
K}$ has eigenvectors relative to the eigenvalue $1$; therefore
\begin{equation*}
e + (3-e)\gamma < 3\gamma +1;
\end{equation*}
thus, if $h$ is so small that
\begin{equation}
\bigl[e +(3-e)\gamma\bigr]C_2 +2\gamma h_1(C_1+LC'_2)\le 
(3\gamma +1)C_2=C'_2,
\label{eq:30}
\end{equation}
$G$ maps $B(0,C'_2)$ into itself; moreover, the Lipschitz constant of
$G$ over this ball is at most equal to $\gamma Lh/2$; if
\begin{equation}
\gamma L h_1 <2,\label{eq:31}
\end{equation}
$G$ is a strict contraction from $B(0,C'_2)$ to itself, which proves
the lemma.
\end{proof}

When $\overline{u}$ belongs to $\partial K$, we need local coordinates 
in which the projection $P_K$ is particularly simple. They are 
defined in the following fashion: we choose a coordinate frame in
$\Er^d$ such that 
\begin{itemize}
\item $\overline{u}=0$;

\item the tangent hyperplane to $\partial K$ at $0$ is the hyperplane
of the first $d-1$ coordinates;

\item the positive direction of the $d$-th coordinate axis points
inside $K$.
\end{itemize}

For a $d$-dimensional vector $x$, we will use the notation
\begin{equation*}
x'=(x_1,\dots, x_{d-1}).
\end{equation*}

Locally, $\partial K$ is a graph over the hyperplane of the
first $d-1$ co\-or\-di\-na\-tes, and it can be parameterized as
\begin{equation*}
\chi(x')=\begin{pmatrix}
x'\\H(x')
\end{pmatrix},
\end{equation*}
where $x'$ belongs to $\Er^{d-1}$, $H$ is of class $C^3$ and $DH(0)$
vanishes. Let $s\mapsto\psi(s,z)$ 
be the parameterization of the geodesic starting at
$z\in \partial K$ with an initial velocity equal to $-N(z)$
which satisfies
\begin{equation}
\vabs[\psi(s,z)]{\frac{\partial \psi}{\partial s}(s,z)} =1.\label{eq:92}
\end{equation}
 Let $\Psi$ be defined by
\begin{equation}
\Psi(x',y)=\psi(-y,\chi(x'));\label{eq:44}
\end{equation}
the function $\Psi$ is of class $C^2$ in a neighborhood of $0$; its
derivative at $0$ has the block representation
\begin{equation}
D\Psi(0,0)=\begin{pmatrix}
\begin{matrix}
\one_{d-1}\\0
\end{matrix} & \Big| & N(0)
\end{pmatrix};\label{eq:93}
\end{equation}
it is invertible, since $N(0)$ does not belong to the tangent plane at
$0$ to $\partial K$. Thus $\Psi$ is a local diffeomorphism from a
neighborhood  $\calu$ of $0$ to a neighborhood $\Psi(\calu)$ of $0$.
In particular, we may assume that $\calu$ contains a compact
neighborhood of $0$ of the form $\overline{\calo} \times [-r_1,r_1]$ where
$\calo$ is an open neighborhood of $0$ in $\Er^{d-1}$. 

The inverse diffeomorphism of $\Psi$ is denoted by $\Phi$, and we
decompose it as
\begin{equation}
\Phi(x)=\begin{pmatrix}
S(x)\\Y(x)
\end{pmatrix},\label{eq:96}
\end{equation}
where $S$ takes its values in $\Er^{d-1}$ and $Y$ takes its values in
$\Er$. If $x$ belongs to 
\begin{equation*}
\calf=\Psi(\overline{\calo} \times[-r_1,r_1]),
\end{equation*}
the projection $P_K(x)$ is given by
\begin{equation}
P_K(x)=\Psi\begin{pmatrix}
S(x)\\Y(x)^+
\end{pmatrix}.\label{eq:97}
\end{equation}
With these preparations, we are able to prove the main local estimates:

\begin{theorem}\label{thr:8}
For all $\overline{u}\in K$, for all $C_2>0$, 
there exist two positive numbers, $r_1<r_2$ and three numbers
 $\tau>0$, 
$h_1>0$ and $C_3<\infty$
such that for all $h\in (0,h_1]$ and all $t_0\in [0,T)$, for
all $U^0$ and $U^1$,  satisfying the
condition $E\bigl(\overline{u},r_1,C_2,h\bigr)$,
$U^m$ is defined in $B(\overline{u}, r_2)$,
for all $m\le \lfloor\tau/h\rfloor$,
and $\vabs{V^m}$ is bounded by  $C_3$ independently of $h$
for $0\le m\le \lfloor\tau/h\rfloor-1$. 
\end{theorem}

\begin{proof}
The theorem decomposes into an easy and a difficult part. The 
easy part is when $\overline{u}$ belongs to the interior of $K$. 

\subsection*{First case: $\overline{u} \in \inter {K}$}

We choose $r_0>0$ as in the proof of lemma~\ref{implicit}:
the ball of center $\overline{u}$ and radius
$2r_0$ is included in $K$. The number $C'_2$ defined by~\eqref{eq:22}
is equal to $4C_2$, and
the numbers $C_1$ and $L$ are given respectively by~\eqref{eq:141}
and~\eqref{eq:23}. 
We choose $r_1=r_0/2$ and $r_2=r_0$.
Assume that $\tau$ satisfies the following 
inequalities:
\begin{equation}
\begin{split}
&0<\tau\bigl(C_1 + LC'_2\bigr)<\min\bigl(C'_2 - C_2,r_0-r_1\bigr),\\
&\tau C_2 
+\frac{\tau^2}{2}\bigl(C_1 + LC'_2\bigr) \le\frac{r_0}{2}.
\end{split}
\label{eq:26}
\end{equation}
Then, if we write 
\begin{equation}
n=\lfloor \tau/h\rfloor,\label{eq:41}
\end{equation}
 we shall prove by
induction that for small enough $h$, there exists a unique solution
of~\eqref{un}, for $0\le m\le n$,
which satisfies the estimate
\begin{equation*}
\vabs{V^m} \le C'_2.
\end{equation*}
We claim that for $h$ small enough, we can find a solution of
\begin{equation}
U^{m+1}-2U^m + U^{m-1} =h^2 F\left(t_m, U^{m-1},U^m,
\bigl(V^m +V^{m-1}\bigr)/2,h\right)
\label{eq:27}
\end{equation}
which satisfies the estimates
\begin{align}
\forall m\in\{0, \dots, n-1\},\quad &\vabs{U^m -\overline{u}}  \le r_1
+ mh\bigl(C_1 + LC'_2\bigr),\label{eq:24}\\
\forall m\in\{0, \dots, n\},\quad &\vabs{V^m} \le C_2 +
mh\bigl(C_1 + LC'_2).\label{eq:25}
\end{align}
In this construction, we seek a solution without considering the
constraints, and we prove eventually that they are satisfied.

It is clear that~\eqref{eq:24} holds for $m\le 2$ and
that~\eqref{eq:25} holds for $m\le 1$. Assume that it holds up to some
exponent $m<n$. Thanks to~\eqref{eq:26}, we have the estimates
\begin{equation*}
\begin{split}
&\vabs{U^{m-2}-\overline{u}}\le r_0, 
\quad\vabs{U^{m-1}-\overline{u}}\le r_0, \quad \vabs{U^m-\overline{u}}\le
r_0,\\
& \vabs{V^{m-2}}\le C'_2,\quad\vabs{V^{m-1}}\le C'_2.
\end{split}
\end{equation*}
Therefore, we may apply lemma~\ref{implicit} with $K=\Er^d$: defining
$C''_2=4C'_2$, 
we can find $h_1$ such that for $0<h\le h_1$, there exists a unique
$U^{m+1}$ such that
\begin{equation*}
\vabs{U^{m+1}-U^m}\le C''_2 h.
\end{equation*}
In particular, if $L'$ is defined by~\eqref{eq:23}, with $C'_2$
replaced by $C''_2$, we infer from~\eqref{eq:27} that
\begin{equation}
\vabs{V^{m}}\le \vabs{V^{m-1}} + h(C_1 + L' C''_2);\label{eq:28}
\end{equation}
therefore, with the help of the induction assumption, we
have the estimate:
\begin{equation*}
\vabs{V^m}\le C_2 + mh(C_1 + LC'_2)+h\bigl(L'C''_2 -L C'_2\bigr).
\end{equation*}
If $h$ satisfies the inequality
\begin{equation*}
C_2 + \tau(C_1 + LC'_2)+h\bigl(L'C''_2 -L C'_2\bigr)\le C'_2,
\end{equation*}
we can see that in fact
\begin{equation*}
\vabs{V^m}\le C'_2,
\end{equation*}
and therefore, instead of~\eqref{eq:28}, we have
\begin{equation*}
\vabs{V^m}\le C_2 + mh\bigl(C_1 + LC'_2\bigr).
\end{equation*}
Therefore, we have also
\begin{equation*}
\vabs{U^{m+1}-\overline{u}}\le \frac{r_0}{2} + (m+1)hC_2
+\frac{m(m+1)h^2}{2}\bigl(C_1 + hC'_2\bigr).
\end{equation*}
Thus,~\eqref{eq:24} and~\eqref{eq:25} hold.
Let us prove that the vector $W^m$ defined by~\eqref{Wn} belongs to
$K$: since
\begin{equation*}
W^m -\overline{u}=U^{m-1} + \frac{2h}{1+e} V^{m-1} +\frac{h^2
F^m}{1+e} -\overline{u} 
\end{equation*}
we have the estimate
\begin{equation*}
\vabs{W^m -\overline{u}}\le r_0 + \frac{2hC'_2}{1+e} +\frac{h^2}{1+e}
\bigl(C_1 + LC'_2); 
\end{equation*}
thus, if $h\le h_1$ and $h_1$ satisfies
\begin{equation*}
\frac{2h_1C'_2}{1+e} +h_1^2 \bigl(C_1 + LC'_2)\le r_0,
\end{equation*}
$W^m$ belongs to $K$, and the sequence $U^m$
satisfies~\eqref{un}. This concludes the proof of the estimates in the
first case. In particular, we can choose $C_3=C'_2$.

\subsection*{Second case: $\overline{u}\in\partial K$}

We define on $\Er^d$ a norm denoted by $\norm{\>\>}$ as follows:
\begin{equation*}
x=\begin{pmatrix}
x'\\x_d
\end{pmatrix},\quad \norm{x}=\max\bigl(\vabs{x'},\vabs{x_d}\bigr).
\end{equation*}
Pick $R_1>0$ such that $\Psi$ is a diffeomorphism from an open neighborhood
\begin{equation*}
\calb=\bigl\{x:\norm{x-\Phi(\overline{u})}\le R_1\bigr\}
\end{equation*}
to its image and such that $\Psi(\calb)$ is included in an euclidean
ball $B(\overline{u},r_0)$ such that $P_K$ is Lipschitz continuous on
$B(\overline{u},2r_0)$; denote by $\gamma$ the Lipschitz constant of
$P_K$ on $B(\overline{u},2r_0)$.

Define $\Lambda$ by
\begin{equation*}
\begin{split}
\Lambda =&\max\biggl\{ 
\sup
\Bigl\{
\frac{\vabs{D\Psi(x)x_1}}{\vnorm{x_1}}:
x\in\calb, x_1\neq 0\Bigr\},\\
&\quad\sup\Bigl\{\frac{\vabs{D\Phi(u)u_1}}{\vnorm{u_1}}:
u\in\Psi(\calb), u_1\neq 0\Bigr\}\biggr\},
\end{split}
\end{equation*}
and
\begin{equation*}
\begin{split}
C_4=&\max\biggl\{\sup\Bigl\{\frac{\vabs{D^2\Psi(x)x_1\otimes
x_2}}{2\vnorm{x_1}\vnorm{x_2}}: 
x\in\calb, x_1\neq 0, x_2\neq 0\Bigr\},\\
&\quad\sup\Bigl\{\frac{\vabs{D^2\Phi(u)u_1\otimes
u_2}}{2\vnorm{u_1}\vnorm{u_2}}: 
u\in\calb, u_1\neq 0, u_2\neq 0\Bigr\}\biggr\}.
\end{split}
\end{equation*}
A continuity argument shows taht the compact set $\Psi(\calb)$
contains the ball of radius $R_1/\Lambda$ about $\overline{u}$.

We will give now a description of the
scheme~\eqref{un},~\eqref{fn} in the new coordinates
$X^m=\Phi(U^m)$. Assume therefore that 
\begin{equation}
U^{m+1}, \> U^m,\>  U^{m-1},\> W^m \text{ and } \frac{U^m + E
U^{m-1}}{1+e} \text{ belong to }\calb.\label{eq:32} 
\end{equation}
We know that~\eqref{un} is equivalent to
\begin{equation}
\frac{U^{m+1} +eU^{m-1}}{1+e}=P_K(W^m);\label{eq:147}
\end{equation}
We map~\eqref{eq:147} by $\Phi$, and we calculate the Taylor expansion
of either side of~\eqref{eq:147} around $U^m$. The left hand side
of~\eqref{eq:147} can be rewritten as
\begin{equation*}
U^m+h\frac{V^m -eV^{m-1}}{1+e}
\end{equation*}
and therefore
\begin{equation}
\begin{split}
&\Phi\left(U^m + h\frac{V^m -eV^{m-1}}{1+e}\right)\\
&\quad=\Phi(U^m)+D\Phi(U^m) h\frac{V^m -eV^{m-1}}{1+e} + I^m\label{eq:148}
\end{split}
\end{equation}
where
\begin{equation*}
\vnorm{I^m}\le C_4\vabs{\frac{h\bigl(V^m -eV^{m-1}\bigr)}{1+e}}^2.
\end{equation*}
But
\begin{equation*}
U^{m+1}=U^m+hV^m,
\end{equation*}
so that another Taylor expansion gives
\begin{equation*}
\Phi\bigl(U^{m+1}\bigr)=\Phi(U^m)+ D\Phi(U^m)hV^m +\hat I^m
\end{equation*}
with
\begin{equation*}
\vnorm{\hat I^m}\le C_4\vabs{hV^m}^2.
\end{equation*}
Thus
\begin{equation}
D\Phi(U^m)hV^m =\Phi(U^{m+1}) -\Phi(U^m) -\hat I^m.\label{eq:149}
\end{equation}
A similar calculation gives
\begin{equation}
\label{eq:150}
-D\Phi(U^m)hV^{m-1} =\Phi(U^{m-1}) -\Phi(U^m) -\tilde I^m,
\end{equation}
with
\begin{equation*}
\vnorm{\tilde I^m}\le C_4\vabs{hV^{m-1}}^2.
\end{equation*}
If we substitute~\eqref{eq:149} and~\eqref{eq:150}
into~\eqref{eq:148}, we find that
\begin{equation*}
\Phi\bigl(U^{m+1}\bigr)=\frac{X^{m+1}+e X^{m-1}}{1+e} -\frac{\overline{I}^m}{1+e}
\end{equation*}
where
\begin{equation*}
\overline{I}^m =\hat I^m +e\tilde I^m -(1+e)I^m,
\end{equation*}
and have the estimate
\begin{equation}
\label{eq:151}
\vnorm{\overline{I}^m}\le C_4 h^2\bigl(\vabs{V^m}^2 + e\vabs{V^{m-1}}^2
+(1+e)^{-1}\vabs{V^m -eV^{m-1}}^2\bigr).
\end{equation}
Consider now the right hand side of~\eqref{eq:147}. By definition of
$V^{m-1}$, we have the identity
\begin{equation}
W^m =U^m +\frac{(1-e)hV^{m-1}+h^2F^m}{1+e},\label{eq:171}
\end{equation}
and a Taylor expansion gives
\begin{equation}
\Phi(W^m)=\Phi(U^m)+D\Phi(U^m)\frac{(1-e)hV^{m-1}+h^2F^m}{1+e} +J^m,\label{eq:152}
\end{equation}
with
\begin{equation*}
\vnorm{J^m}\le C_4\vabs{\frac{(1-e)hV^{m-1}+h^2F^m}{1+e}}^2.
\end{equation*}
We substitute~\eqref{eq:150} into~\eqref{eq:152}, and we obtain
\begin{equation*}
\Phi(W^m)=\frac{2X^m -(1-e)X^{m-1} +h^2D\Phi(U^m)F^m}{1+e}
+\frac{\overline{J}^m}{1+e}, 
\end{equation*}
where
\begin{equation*}
\overline{J}^m =(1+e)J^m +(1-e)\tilde I^m,
\end{equation*}
so that
\begin{equation*}
\vnorm{\overline{J}^m}\le C_4\left[\frac{\vabs{(1-e)hV^{m-1}
+h^2F^m}^2}{1+e}+(1-e)h^2\vabs{V^{m-1}}^2\right]. 
\end{equation*}
We have to estimate $\vnorm{\overline{I}^m}+\vnorm{\overline{J}^m}$;
by elementary inequalities, 
\begin{equation*}
\begin{split}
&\vnorm{I^m}+\vnorm{J^m} \le C_4h^2\Bigl[\frac{2(1-e)^2\vabs{V^{m-1}}^2
+2h^2\vabs{F^m}^2}{1+e} \\&\quad+(1-e)\vabs{V^{m-1}}^2 +
\vabs{V^m}^2+e\vabs{V^{m-1}}^2 +\frac{2\vabs{V^m}^2
+2e^2\vabs{V^{m-1}}^2}{1+e}\Bigr].
\end{split}
\end{equation*}
The coefficient of $\vabs{V^{m-1}}^2$ in the above bracket is
\begin{equation*}
\frac{2(1-e)^2}{1+e}+1 +\frac{2e^2}{1+e},
\end{equation*}
and since for $e\in [0,1]$, $(1-e)^2\le 1-e^2$, this coefficient is at
most equal to $3$. The coefficient of $\vabs{V^{m}}^2$ in the same
bracket is at most equal to 
\begin{equation*}
1+\frac{2}{1+e},
\end{equation*}
which is also at most equal to $3$. Therefore
\begin{equation}
\label{eq:153}
\vnorm{\overline{I}^m}+\vnorm{\overline{J}^m} \le C_4 h^2\bigl[3\vabs{V^m}^2
+3\vabs{V^{m-1}}^2 + 2h^2\vabs{F^m}^2\bigr].
\end{equation}
Thanks to the properties of $P_K$,
\begin{equation}
\Phi\left(\frac{U^{m+1}+eU^{m-1}}{1+e}\right)=\Phi(P_KW^m)=\begin{pmatrix}
S(W^m)\\Y(W^m)^+
\end{pmatrix}.\label{eq:180}
\end{equation}
Define
\begin{equation*}
s^m=S(U^m)=\bigl[X^m\bigr]', \quad y^m =Y(U^m)=X^m_d.
\end{equation*}
In these new coordinates, we have
\begin{align}
\label{eq:154} s^{m+1} -2s^m +s^{m-1}=h^2\kappa^m,\\
\label{eq:155}y^{m+1}+ey^{m-1}-\bigl(2y^m -(1-e)y^{m-1}\bigr)^+=h^2\lambda^m,
\end{align}
and $\kappa^m$ and $\lambda^m$ are given by
\begin{equation*}
\begin{split}
h^2\kappa^m &=\bigl[h^2D\Phi(U^m)F^m +\overline{I}^m +\overline{J}^m\bigr]',\\
h^2\lambda^m&=\bigl[2y^m -(1-e)y^{m-1} +\bigl(h^2 D\Phi(U^m)F^m
+\overline{J}^m\bigr)_d\bigr]^+ \\ 
&\qquad\qquad-\bigl(2y^m -(1-e)y^{m-1}\bigr)^+ + \overline{I}^m_d.
\end{split}
\end{equation*}
Therefore, we have the estimates:
\begin{equation}
\max\bigl(\vabs{\kappa^m},\vabs{\lambda^m}\bigr)\le \Lambda\vabs{F^m}
+C_4\bigl(3\vabs{V^m}^2 +3\vabs{V^{m-1}}^2 + 2h^2\vabs{F^m}^2\bigr).\label{eq:45}
\end{equation}
We define $\xi^m$ and $\zeta^m$ by
\begin{equation*}
\xi^m=\begin{pmatrix}
\sigma^m\\\eta^m
\end{pmatrix}=\frac{X^{m+1}-X^m}{h},\quad \zeta^m=\begin{pmatrix}
\kappa^m\\\lambda^m
\end{pmatrix}.
\end{equation*}
Let now $q$ be a number which satisfies
\begin{equation*}
q > \Lambda C_2.
\end{equation*}
Let 
\begin{equation}
C'_2=\Lambda q(3\gamma +1),\label{eq:50}
\end{equation}
and let $C_1$
and $L$ be respectively as in~\eqref{eq:141} and~\eqref{eq:23}.
If we assume beyond~\eqref{eq:32} that
\begin{equation}
\max\bigl(\vabs{V^{m-1}},\vabs{V^m}\bigr)\le C'_2,\label{eq:33}
\end{equation}
we have the estimate
\begin{equation*}
\vabs{F^m}\le C_1 + L\bigl(\vabs{V^m}+\vabs{V^{m-1}}\bigr)/2;
\end{equation*}
by elementary inequalities,
\begin{equation*}
\vabs{F^m}^2 \le 2C_1^2 + L^2\bigl(\vabs{V^m}^2+\vabs{V^{m-1}}^2\bigr),
\end{equation*}
and therefore, if we define
\begin{equation*}
C_5=\frac{\Lambda}{2} +C_1^2\bigl(\Lambda +4h_1^2 C_4\bigr),\quad
C_6=\Bigl(\Bigl(\frac{\Lambda}{2} +2h_1^2 C_4\Bigr)L^2 + 3C_4\Bigr)\Lambda^2,
\end{equation*}
we have shown that under assumptions~\eqref{eq:32} and~\eqref{eq:33},
the following inequality holds:
\begin{equation}
\vnorm{\zeta^m} \le C_5 + C_6\bigl(\vnorm{\xi^m}^2
+\vnorm{\xi^{m-1}}^2\bigr).\label{eq:37} 
\end{equation}
Let $\tau$ be a number which satisfies the following inequalities:
\begin{equation}
\begin{split}
&\tau>0, \quad \Lambda C_2 +(2\tau C_5 + 2C_6 q^2)<q,\\
&0<\rho= \frac{R_1}{2} -\tau \Lambda C_2 - 2\tau^2\bigl(C_5 + 2 C_6q^2\bigr).
\end{split}
\label{eq:36} 
\end{equation}
Assume that initially
\begin{equation}
\max_{j=0,1,2}\vabs{U^j-\overline{u}}\le R_1/(2\Lambda), \quad 
\max_{j=0,1}\vabs{U^{j+1}-U^j}\le C_2h.\label{eq:34}
\end{equation}
We will prove by induction that if $n=\lfloor \tau/h\rfloor$, then for
all $m\le n$
\begin{equation}
\begin{split}
&\forall l\in\{0, \dots, m\},\quad \vnorm{X^l -\Phi(\overline{u})} \le
\vnorm{X^0-\Phi(\overline{u})} \\&\qquad+ lh \Lambda C_2 +2l(l-1)h^2\bigl(C_5
+2C_6 q^2\bigr),\\
&\forall l\in\{0, \dots, m-1\}, \quad \vnorm{\xi^l} \le \Lambda C_2 +
2lh\bigl(C_5 +2 C_6 q^2\bigr).
\end{split}\label{eq:35}
\end{equation}
For $m\le 2$, assumptions~\eqref{eq:34} guarantee that~\eqref{eq:35}
holds. 
The induction hypotheses imply that
\begin{equation*}
\max_{j=m-2,m-1,m} \vnorm{X^j -\Phi(\overline{u})} \le R_1-\rho, \quad
\max_{j=m-2,m-1} \vnorm{X^{j+1}-X^j}\le qh,
\end{equation*}
and therefore, $U^{m-2}$, $U^{m-1}$ and $U^m$ belong to
$\Psi(\calb)$. We may apply lemma~\eqref{implicit} which guarantees
the existence of $U^{m+1}$ such that
\begin{equation}
\vabs{U^{m+1}-U^m} \le (3\gamma +1)\Lambda qh.\label{eq:39}
\end{equation}
The ball of radius $\rho/\Lambda$ about $U^m$ is included in
$\Psi(\calb)$; thus, if
\begin{equation*}
(3\gamma +1)\Lambda^2 qh \le \rho,
\end{equation*}
$U^{m+1}$ also belongs to $\Psi(\calb)$. Similarly,
\begin{equation*}
\vabs{\frac{U^{m+1}+eU^{m-1}}{1+e} - U^m} \le (3\gamma +1)\Lambda qh,
\end{equation*}
and if
\begin{equation*}
(3\gamma +1°\Lambda^2 qh\le \rho,
\end{equation*}
$(U^{m+1}+ eU^{m-1})/(1+e)$ belongs to $\Psi(\calb)$.
Finally, thanks to the definition~\eqref{eq:50} of $C'_2$, we have the
inequality:
\begin{equation*}
\vabs{\frac{U^{m+1}-U^{m-1}}{2h}}\le \frac{(3\gamma +2)\Lambda q}{2}
\le C'_2;
\end{equation*}
in virtue of the definitions~\eqref{eq:23} of $L$ and~\eqref{eq:141}
of $C_1$, we have 
\begin{equation*}
\vabs{W^m - U^m} \le \frac{1-e}{1+e} \Lambda q h + \frac{h^2}{1+e}
\bigl(LC'_2 + C_1\bigr).
\end{equation*}
Once again, if
\begin{equation*}
\frac{1-e}{1+e} \Lambda q h + \frac{h^2}{1+e}
\bigl(LC'_2 + C_1\bigr)\le \frac{\rho}{\Lambda},
\end{equation*}
$W^m$ belongs to $\Psi(\calb)$. Thus~\eqref{eq:32} holds and we may
apply the argument that followed.

By definition of $\sigma^m$ and $\eta^m$, we have the inequalities
\begin{align*}
\vabs{\sigma^m}&\le \vabs{\sigma^{m-1}}+h\vnorm{\zeta^m},\\
\intertext{and thanks to lemma~\ref{monodim}}
\vabs{\eta^m}&\le \max\bigl(\vabs{\eta^{m-1}},e\vabs{\eta^{m-2}}\bigr)
+ h\vnorm{\zeta^m} +h\vnorm{\zeta^{m-1}};
\end{align*}
hence we infer that
\begin{equation}
\vnorm{\xi^m} \le \max\bigl(\vnorm{\xi^{m-1}},
e\vnorm{\xi^{m-2}}\bigr)
+ h\vnorm{\zeta^m} + h\vnorm{\zeta^{m-1}},\label{eq:51}
\end{equation}
and thanks to the induction hypothesis
\begin{equation}
\begin{split}
\vnorm{\xi^m} &\le \Lambda C_2 + 2mh(C_5 + 2C_6 q^2) -hC_6q^2 \\
&\quad+ h(C_5 +
C_6q^2) + hC_6 \vnorm{\xi^m}^2.
\end{split}
\label{eq:40}
\end{equation}
The equation in $a$
\begin{equation} 
hC_6a^2 - a +2mh(C_5 + 2C_6 q^2) -hC_6q^2 + h(C_5 + C_6q^2) +\Lambda
C_2=0\label{eq:38} 
\end{equation}
has two distinct real roots if
\begin{equation*}
\Delta = 1-4hC_6 \bigl(2mh(C_5 + 2C_6 q^2) -hC_6q^2 + h(C_5 +
C_6q^2)+\Lambda C_2\bigr)
\end{equation*}
is strictly positive; but this is always true if $0<h\le h_1$ and
\begin{equation*}
1>4h_1 C_6(2\tau \bigl(C_5+2C_6 q^2) +\Lambda C_2\bigr).
\end{equation*}
The smallest of the two roots of~\eqref{eq:38} is inferior to $q$,
since the substitution $a=q$ in~\eqref{eq:38} gives a negative left
hand side; the largest of these two roots is at least equal to
$1/(2hC_6)$; but relation~\eqref{eq:39} implies
\begin{equation*}
\vnorm{\xi^m}\le \Lambda C'_2;
\end{equation*}
thus, if
\begin{equation*}
h_1\le \frac{1}{2C_6\Lambda C'_2},
\end{equation*}
relation~\eqref{eq:40} implies $\vnorm{\xi^m}\le q$; if we substitute
this inequality in the right hand side of~\eqref{eq:40}, we find that
the second inequality in~\eqref{eq:35} holds for $l=m$; the first
inequality in~\eqref{eq:35} for $l=m+1$ holds immediately, and the
induction is proved. Thus, we can take as an upper bound of
$\vabs{V^m}$ the number $C_3=\Lambda q$; we can take also
$r_1=R_1/(2\Lambda)$ and $r_2=\Lambda R_1$.
\end{proof}

If we put together theorem~\ref{thr:8} and lemma~\ref{implicit0}, we
obtain an existence result:

\begin{theorem}\label{thr:12}For all $(t_0, u_0,M(u_0)v_0)\in \mathbb{D}$,
for all $U^1$ satisfying~\eqref{u1}, there exists $\tau>0$,
$C_3<\infty$ and $h_1$ 
such that for all $h\in(0,h_1]$, there exists a unique solution of~\eqref{un}
and~\eqref{fn} for all $m\le \lfloor\tau/h\rfloor-1$, which satisfies
the estimate
\begin{equation}
\forall l\le n-1, \quad\vabs{V^l}\le C_3.\label{eq:42}
\end{equation}
\end{theorem}

\begin{proof} Let us check that $U^0$ and $U^1$ satisfy
condition $E(u_0,r_1,C_2,h)$. Lemma~\ref{implicit0} and
assumption~\eqref{u1} on $U^1$ imply that 
\begin{equation*}
\vabs[u_0]{U^1 -u_0} \le h \bigl(\vabs[u_0]{z(h)} + \vabs[u_0]{v_0}\bigr)
\end{equation*}
and
\begin{equation*}
\vabs[u_0]{U^2 -u_0} \le h\bigl(2\vabs[u_0]{v_0} + 1\bigr).
\end{equation*}
Choose $C_2\ge \bigl(4\vabs[u_0]{v_0} + 1\bigr)\vnorm{M(u_0)}$; $U^0$
and $U^1$ satisfy condition $E(u_0,r_1,C_2,h)$ for small enough values
of $h$. 
Then, it is clear that theorem~\ref{thr:8} applies.
\end{proof}

It is convenient to give a uniformized version of theorem~\ref{thr:8}:

\begin{theorem}\label{thr:9} For all compact subset $\calc$ of $K$, for all
$C_2>0$, there exist positive numbers $r_1$, $r_2>r_1$, $\tau$, $C_3$, and $h_1$ such
that for all $t_0\in[0,T)$, for all $\overline{u}\in\calc$, for all $h\le
h_1$ and for all $U^0$
and $U^1$ satisfying condition $E(\overline{u},r_1,C_2,h)$
relations~\eqref{un} and~\eqref{fn} define uniquely 
under condition~\eqref{eq:42} the vectors $U^m$
for $2\le m\le \lfloor\min(\tau,T-t_0)/h\rfloor]$.
\end{theorem}

\begin{proof}
Any element $u$ of $\calc$ is included in an open ball
$\inter{B(u,r_1(u))}$ such that theorem~\ref{thr:8} holds.
We cover $\calc$ by a finite number of balls
$\inter{B(u_j,r_1(u_j)/2)}$ with associated numbers $r_2(u_j)$,
$\tau(u_j)$, $h_1(u_j)$ and $C_3(u_j)$. If we let 
\begin{equation*}
r_1=\frac{1}{2} \min\{r_1(u_j):1\le j \le J\};
\end{equation*}
then any $\overline{u}\in \calc$ belongs to a ball $B(u_j,r_1(u_j)/2)$, and in
particular, $B(\overline{u},r_1)$ is included in $B(u_j,r_1(u_j))$. 
If we take
\begin{equation*}
{\tau}=\min_j \tau(u_j), \> {r_2}=\max_j r_2(u_j)
\quad {h}_1=\min_j h_1(u_j), \> {C_3}=\max_j C_3(u_j),
\end{equation*}
it is immediate that the theorem holds, thanks to theorem~\ref{thr:8}.
\end{proof}

\section{Estimates on the acceleration\label{sec:Estim-accel}}

In this section and the three following ones, we assume that there
exist strictly positive numbers $\tau$, $C_3$ and $h_1$, and a
subsequence of times steps to which correspond
solutions of the numerical scheme defined
by~\eqref{u0}, \eqref{u1}, \eqref{un} and~\eqref{fn}, which
satisfy the estimate, for all $h\le h_1$:
\begin{equation}
\forall l\in \{0,P -1\},\quad \vabs{U^{l+1}-U^l}\le
C_3 h\label{eq:168} 
\end{equation}
where
\begin{equation*}
P= \lfloor \tau/h\rfloor
\end{equation*}
Here we estimate the discrete total variation of the sequence
$\bigl(V^m\bigr)_m$. 
It is also convenient to define the function $w_h(t)$ on $[t_0,
t_0+\tau]$ by
\begin{equation*}
\begin{split}
&w_h(t_m)=W^m, \quad\text{$w_h$ is continuous and it is affine}\\
&\text{on each
interval $[t_m, t_{m+1})$, and constant on $[t_P,t_0+\tau]$.}
\end{split}
\end{equation*}

\begin{theorem}\label{thr:10}Under assumption~\eqref{eq:168}, there exists a
constant $C_7$ such that for all $h\le h_1$: 
\begin{equation} 
\sum_{m=1}^{P-1}\vabs{V^m-V^{m-1}}\le C_7.\label{eq:181}
\end{equation} 
\end{theorem}

\begin{proof} Let $\calc$ be the compact set $K\cap
B(u_0,C_3\tau)$ and let
$r_1$ be as in theorem~\ref{thr:9}; cover $\calc$ with a finite number of
balls $B(u_j,r_1/4)$; observe that, thanks to
Ascoli--Arzel\'a's theorem, the set $\calw$ of functions $(w_h)_{0<h_1\le h}$
is relatively compact in $C^0([t_0,t_0+\tau])$. The
set of limit points of $(w_h)_{0<h\le h_1}$ as $h$ tends to $0$ is
also a compact set, which we shall denote by $\calw_{\infty}$. There
exists a finite subset $w^1, \dots, w^I$ of $\calw_\infty$ such that 
\begin{equation*}
\forall w\in \calw_\infty: \inf\{\norm[{C^0[t_0, t_0+\tau]}]{w
-w^i}: 1\le i \le I\} \le r_1/4.
\end{equation*}
For each $i\in \{1,\dots,I\}$, it is possible to find a
finite increasing sequence of times
\begin{equation*}
0=\tau(i,0)<\dots <\tau(i,k)<\dots <\tau(i,\kappa(i))=\tau
\end{equation*}
such that 
\begin{equation*}
w^i([\tau(i,k),\tau(i,k+1)])\subset B\bigl(u_{j(i,k)},
 r_1/4\bigr). 
 \end{equation*}
Thus, for all $w\in\calw_\infty$, 
\begin{equation*}
w([\tau(i,k),\tau(i,k+1)])\subset B\bigl(u_{j(i,k)},
 r_1/2\bigr).
\end{equation*}
Therefore, we can
decrease $h_1$ so that  
\begin{equation*}
\begin{split}
&\forall h\in(0,h_1], \quad\exists i\in\{1,\dots, I\},\quad\forall
k\in\{1,\dots ,\kappa(i)\}, \\&\forall t\in [\tau(i,k),\tau(i,k+1)]\quad
w_h(t)\in
B\bigl(u_{j(i,k)}, 3r_1/4\bigr),
\end{split}
\end{equation*}
and thanks to~\eqref{eq:171} and to~\eqref{eq:168}, we can
decrease $h_1$ such that
\begin{equation*}
\begin{split}
&\forall h\in(0,h_1], \quad\exists i\in\{1,\dots, I\},\quad\forall
k\in\{1,\dots ,\kappa(i)-1\},\\
&\qquad \forall l\in\{\lfloor
\tau(i,k)/h\rfloor, \dots \lfloor \tau(i,k+1)/h\rfloor\}, \quad U^l
\in B\bigl(u_{j(i,k)}, r_1\bigr).
\end{split}
\end{equation*}
We simplify the notations by letting
\begin{equation*}
P=\lfloor \tau(i,k)/h\rfloor,\quad Q= \lfloor \tau(i,k+1)/h\rfloor,
\end{equation*}
and we take $C_1$ be as in~\eqref{eq:141}, where
$\overline{u}$ is set equal to $u_{j(i,k)}$, $r_0$ is set
equal to  $r_1$ and $C'_2$ is set equal to $C_3$.

Now, we have to consider two cases: 

\subsection*{First case: $B\bigl(\overline{u},
r_1\bigr) \cap \partial K=\emptyset$}
We have the inequality
\begin{equation}
\vabs{F^m}\le C_1+ LC_3,\label{eq:178}
\end{equation}
hence, thanks to~\eqref{un}, we have the inequality
\begin{equation*}
\vabs{V^m -V^{m-1}}\le h\bigl(C_1+ LC_3\bigr),
\end{equation*}
and therefore
\begin{equation}
\sum_{m=P+1}^{Q}\vabs{V^m -V^{m-1}} \le
\bigl(\tau(i,k+1)+2h-\tau(i,k)\bigr)\bigl(C_1+ LC_3\bigr). \label{eq:176}
\end{equation}

\subsection*{Case 2: $B\bigl(\overline{u},
r_1\bigr) \cap \partial K\neq\emptyset$}
We observe that thanks to~\eqref{eq:37}, we have the estimate
\begin{equation}
\forall m \in \{P+1,\dots, Q-1\}, \quad
\max\bigl(\vabs{\kappa^m},\vabs{\lambda^m}\bigr) \le C_{9},\label{eq:177}
\end{equation}
where
\begin{equation*}
C_{9} =C_5 + 2 C_6 \Lambda ^2C_3^2.
\end{equation*}

The estimates on the first $d-1$ components of the velocity in the
straightened coordinates are immediate:
\begin{equation}
\sum_{m=P+1}^Q \vabs{\frac{s^{m+1}-s^{m}}{h} - \frac{s^m -s^{m-1}}{h}}
\le \bigl(\tau(i,k+1)+2h-\tau(i,k)\bigr) C_{9}.\label{eq:174}
\end{equation}
In order to estimate the last coordinate, we partition $\{P+1,\dots,
Q\}$ as follows: 
\begin{equation*}
\begin{split}
\calp&=\bigl\{m\in\{P+1,\dots,Q\}:2y^m -(1-e)y^{m-1}<0\}, \\
\calp'&=\{P+1,\dots,Q\}\setminus \calp.
\end{split}
\end{equation*}
We write $\calp$ as an union of discrete intervals:
\begin{equation*}
\calp=\bigcup_{l=1}^\ell \{p(l),\dots,q(l)\}, \quad p(l)-1\notin \calp,
\quad q(l)+1\notin\calp. 
\end{equation*}
If $\eta^m$ is defined as in~\eqref{eq:107}, we observe that for
$m\in\calp'$,
\begin{equation*}
\vabs{\eta^m -\eta^{m-1}}\le C_{9} h,
\end{equation*}
so that
\begin{equation*}
\sum_{m\in\calp'} \vabs{\eta^m -\eta^{m-1}} \le h C_{9} \vabs{\calp'}.
\end{equation*}
If $m$ belongs to $\calp$, we observe that
\begin{equation}
\eta^m -\eta^{m-1} =h\lambda^m +\bigl(2y^m -(1-e)y^{m-1}\bigr)^-\label{eq:170}
\end{equation}
and therefore, by the triangle inequality,
\begin{equation*}
\vabs{\eta^m -\eta^{m-1}}\le h C_{9}+ \bigl(2y^m -(1-e)y^{m-1}\bigr)^-,
\end{equation*}
and using~\eqref{eq:170} again,
\begin{equation}
\vabs{\eta^m -\eta^{m-1}} \le 2h C_{9} +\eta^m - \eta^{m-1}.\label{eq:172}
\end{equation}
We observe that we have the elements of a telescoping sum: we
sum \eqref{eq:172} for $m$ varying from $p(l)+1$ to $q(l)$, and we
obtain
\begin{equation}
\sum_{m=p(l)}^{q(l)} \vabs{\eta^m -\eta^{m-1}} \le 2h C_{9}
\bigl(q(l)-p(l)\bigr) + \eta^{q(l)}-\eta^{p(l)}.\label{eq:173}
\end{equation}
Now, we sum~\eqref{eq:173} from $l=1$ to $\ell$, which yields
\begin{equation*}
\begin{split}
&\sum_{l=1}^\ell \sum_{m=p(l)}^{q(l)} \vabs{\eta^m -\eta^{m-1}}
\\& \quad \le 
2 C_{9} (Q-P) - \eta^{p(1)}+ \eta^{q(\ell)} - \sum_{l=2}^\ell
\eta^{p(l)} -\eta^{q(l-1)}.
\end{split}
\end{equation*}
But the terms $\eta^{p(l)} -\eta^{q(l-1)}$ can be estimated, since
they correspond to a summation over $\calp'$:
\begin{equation*}
\vabs{\eta^{p(l)} -\eta^{q(l-1)}} \le C_{9} h \bigl(p(l) - q(l)\bigr).
\end{equation*}
Therefore, we have proved that
\begin{equation*}
\sum_{m=P+1}^Q \vabs{\eta^m -\eta^{m-1}} \le 3 C_{9} h(Q-P) + 2 C_3
\Lambda. 
\end{equation*}
Summarizing this relation with~\eqref{eq:174}, we can see that
\begin{equation}
\sum_{m=P+1}^Q \vabs{V^m - V^{m-1}} \le \Lambda C_{9}\bigl(4
\tau(i,k+1) -4\tau(i,k)+8h \bigr)+ 2C_3\Lambda.\label{eq:175}
\end{equation}
Relations~\eqref{eq:176} and~\eqref{eq:175} do not depend on $h\le
h_1$; since we have only a finite number of these estimates, the
theorem is proved.
\end{proof}

\section{Variational properties of the limit of the numerical scheme}
\label{sec:Variational}

In this section, we work under the assumption~\eqref{eq:168}.
Recall that $n=\lfloor \tau/h\rfloor$.
We define a function $u_h$ by affine interpolation, as follows:
\begin{equation}
\begin{cases}
u_h (t) =
 U^m + (t-mh) \frac{\displaystyle U^{m+1} -U^m }{\displaystyle h} 
\\
\qquad\qquad \text{ for $t
 \in [mh, (m+1)h \bigl)$, $0 \le m \le n-1$,} 
\\
u_h(t) = U^{n} \quad \text{ for $t \in [nh,\tau]$.}
\end{cases}\label{eq:79}
\end{equation}

We also define a measure $F_h$ as the following sum of Dirac masses:
\begin{equation}
F_h(t)=\sum_{m=1}^n h F^m \delta(t-mh).\label{eq:80}
\end{equation}

In this section we prove that the sequence $(u_h)_{h}$ converges in an
appropriate sense to a function $u$ which satisfies~\eqref{eq:2}
to~\eqref{eq:8} with $\tau$ instead of $T$. 
We delay the proof of~\eqref{eq:6}, the
transmission condition at impacts, to a later section.

There are three steps in the convergence proof: the first is to prove
that the limit $u$ exists in an appropriate sense and takes its values
in $K$; in the second step, we show that $\dot u_h$ is of bounded
variation uniformly in $h$ and that $F_h$ converges to
$M(u)^{-1}f(\cdot, u, M(u)\dot u)$ weakly in the space of
$\Er^d$-valued measures.
The last step is the characterization of the measure $\mu=M(u)\ddot
u-f(\cdot, u, M(u)\dot u)$: there we show that $\mu$ satisfies
conditions~\eqref{eq:3}, \eqref{eq:4} and~\eqref{eq:5}.

\begin{lemma}\label{l4.8} From all sequence of functions $(u_h)_h$
indexed by a sequence $h$ tending to $0$, it is possible to extract a
subsequence, still denoted by $(u_h)_h$ such that
\begin{align}
u_h &\to u \quad \hbox{\rm in $C^0 ([t_0, t_0+\tau])$ strong},\label{eq:76}\\
\dot u_h &\to \dot  u \quad \text{\rm in $L^{\infty} ([t_0,
t_0+\tau])$ weak  
*}.\label{eq:77}
\end{align}
The function $u$ takes its values in $K$.
\end{lemma}

\begin{proof}
Thanks to assumption~\eqref{eq:168}, we know
that $(u_h)_{0<h \le  h_1}$ is uniformly Lipschitz continuous over
$[t_0, t_0+\tau]$. Therefore, we may extract a subsequence,
still
denoted by $u_h$, such that \eqref{eq:76} and~\eqref{eq:77} hold.
Thus $u$ belongs to $W^{1,\infty}([t_0,t_0+\tau]) \cap
C^0([t_0,t_0+\tau])$,
which means that $u$ 
is a Lipschitz continuous function \cite{brezisfonc}.
For all $m$ belonging to $ \{1, \ldots, n \}$, we have:
\begin{equation}
Z^m = \frac{U^{m+1} + e U^{m-1} }{ 1+e} = U^m + h \frac{ V^m - e
V^{m-1}}{ 1+e}, 
\end{equation}
hence 
$U^m = Z^m - h(V^m - e V^{m-1})/(1+e)$. 
By definition of the scheme, we have $Z^m=P_K(W^m)$(\eqref{eq:78}),
and thus $Z^m$ belongs to $K$. It follows that, for all $m\in \{1,
\dots, n\}$, the euclidean distance between $U^m$ and $K$ can be
estimated as follows:
\begin{equation}
\min\bigl\{\vabs{U^m -u}: u\in K\bigr\} \le h\vabs{V^m - e
V^{m-1}}/(1+e)\le hC_3.\label{eq:169}
\end{equation}
Thanks
to the definition~\eqref{eq:79}, we can see that for all $t\in
[t_0,t_0+\tau]$ 
the euclidean distance between $u_h(t)$ and $K$ is estimated by
$2hC_3$. This allows us to pass to the limit when $h$ tends to $0$ and
to conclude.
\end{proof}

Next lemma describes the convergence of the measures involved in our
problem; we denote by $M^1\bigl((t_0, t_0+\tau) \bigr)$ the space of
bounded measures over $(t_0,t_0+\tau)$ with values in $\Er^d$.

\begin{lemma}The measures $\ddot u_h$ and $F_h$ converge weakly in
$M^1\bigl((t_0,t_0+\tau) \bigr)$ respectively to $\ddot u$ and
$M(u)^{-1}f(\cdot, u,M(u)\dot u)$.\label{thr:2}
\end{lemma}

\begin{proof}
The measure $\ddot u_h$ is a sum of Dirac measures on 
$(t_0,t_0+\tau)$,
more precisely, we have:
\begin{equation}
\ddot u_h (t) =  \sum _{m=1}^{n} (V^m - V^{m-1}) \delta(t-mh) -V^{n}
 \delta (t -  (n+1)h ),
\end{equation}
and the total variation of $\dot u_h$ on $(t_0,t_0+\tau)$ is estimated by
\begin{equation}
TV( \dot u_h ) \le \sum _{m=1}^{n} \vabs{V^m - V^{m-1} } + \vabs{V^{
n}}.
\end{equation}
Theorem~\ref{thr:10} implies that $(\dot u_h)_{0< h
\le h_1}$ is a bounded family in $BV\bigl((t_0,t_0+\tau)\bigr)$,
the space of functions 
of bounded variation over $(t_0,t_0+\tau)$, with values in $\Er^d$. 
Using Helly's 
theorem, we can extract another subsequence $\bigl(\dot u_h\bigr)_h$
which converges, except perhaps on a countable set of points, to a 
function of  bounded variation. Hence
\begin{equation*}
\dot u \in BV\bigl((t_0,t_0+\tau)\bigr).
\end{equation*}
Moreover,
\begin{equation*}
\ddot u_h \to \ddot u \quad \text{ weakly in  $M^1\bigl((t_0,t_0+\tau)
\bigr)$.} 
\end{equation*}
 
Le\-bes\-gue's theorem implies that 
$\dot u_h$ converges to $\dot u$ in $L^1 \bigl(t_0,t_0+\tau)$. We extend
$\dot u_h$ and $\dot u$ to $\Er$ by $0$ outside of $(t_0,t_0+\tau)$
and still denote the respective extensions by $\dot u_h$ and $\dot u$.
The set $\{\dot u_h: h\in (0, h_1]\}\cup \{\dot u\}$ is a compact
subset of $L^1(\Er)$. 
The classical characterization of compact subsets of $L^1(\Er)$ \cite{dunfordschwartz2}
implies that
\begin{equation}
\lim_{\mbox{{}}\theta \to 0} \sup_{0 < h \le h_1} \int_{\Er} \vabs{\dot u_h
(t- \theta ) - \dot u_h(t) 
} \,dt =0.
\end{equation}
Letting $\theta = h$, we can see that $\dot u_h (\cdot -h)$ converges to
$\dot u$  in $L^1 \bigl( 
\Er )$.  
Let us define an approximate velocity $v_h$ on $\Er$ by
\begin{equation}
v_h (t) =  \frac{ \dot u_h (t-h +0) + \dot u_h ( t+0)}{ 2}.
\end{equation}
The sequence $v_h$ converges to $\dot u$ in $L^1 \bigl( \Er
\bigr)$. Moreover, for all $t \in [t_m, t_{m+1})$ and for all $m
\in\{1,\ldots, n\}$,  
we have the identity
\begin{equation}
v_h(t) = \frac{ V^m + V^{m-1}}{ 2}.
\end{equation}
We have immediately the following estimates for all $t \in
(t_0,t_0+\tau)$ and all $h\in (0,h_1]$: 
\begin{equation}
\vabs{v_h(t)} \le C_3 ;\quad\vabs{u_h(t)-u_0}\le C_3(t-t_0)\le
C_3\tau.\label{eq:83}
\end{equation}
Let $\psi$ be a continuous function over $[0,T]$ with compact support
included in $(t_0,t_0+\tau)$. 
For all small enough $h$, the support of $\psi$ is included in
$[t_0+h,t_0+nh]$. 
The duality product $\langle F_h, \psi \rangle$
has the expression
\begin{equation}
\langle F_h, \psi \rangle = \sum_{m=1}^{n} h\psi(t_0+mh)^T F^m.\label{3.28}
\end{equation}
We wish to compare the expression~\eqref{3.28} to
\begin{equation}
\int_{t_0}^{t_0+\tau} \psi^T M(u)^{-1}f(\cdot, u,M(u)\dot u)\,
dt.\label{eq:81} 
\end{equation}
We compare the right hand side of~\eqref{3.28} which is basically a numerical
quadrature by the formula of rectangles to an appropriate integral.
Let us rewrite the individual terms of the right hand side of
\eqref{3.28} as 
\begin{equation}
h \psi (t_m)^T F^m = \int_{t_m}^{t_{m+1}} \psi (t)^T F^m  \, dt + 
\int_{t_m}^{t_{m+1}}  \bigl(\psi (t_m) -\psi (t)\bigr)^T  F^m  \,
dt.\label{3.29}
\end{equation}
Consider the second term on the right hand side of~\eqref{3.29}:
we have already proved (see relations \eqref{eq:178} and \eqref{eq:177})
that there exists a constant $C_{8}$
independent of $m$ and $h\le h_1$ such that
\begin{equation}
\max_{0\le m\le n}\vabs{F^m} \le C_{8}.\label{eq:75}
\end{equation}
Denoting by 
$\omega_\psi$ the modulus of continuity of $\psi$ 
we can see that
\begin{equation}
\vabs{\int_{t_m}^{t_{m+1}}   \bigl( \psi (t_m) -\psi (t) \bigr)^T F^m
\, dt} \le C_{8} \omega_{\psi} (h) h.\label{eq:85}
\end{equation}
We consider now the first term on the right hand side of~\eqref{3.29},
which we would like to compare to expression~\eqref{eq:81}. 
Thanks to the consistance assumption~\eqref{consistance} have the
following 
inequalities, for all $t\in [t_m, t_{m+1})$, and all $n\in
\{1,\dots, n\}$:
\begin{equation*}
\begin{split}
&\vabs{F^m- M(u_h(t))^{-1} f\bigl(t,u_h(t), M(u_h(t))v_h(t)\bigr)}\\
&\qquad\le \vabs{F\bigl(t_m,U^m, U^{m-1},v_h(t_m),h\bigr)- F\bigl(t_m,u_h(t),
u_h(t),v_h(t_m),h\bigr)} \\
&\qquad+ \vabs{F\bigl(t_m,u_h(t),
u_h(t),v_h(t_m),h\bigr)-F\bigl(t_m,u_h(t), u_h(t),v_h(t),0\bigr)}\\
&\qquad+\bigl|M(u_h(t))^{-1} \bigl[f\bigl(t_m,u_h(t),M(u_h(t))
v_h(t)\bigr) -\\
&\qquad\qquad
f\bigl(t,u_h(t),M(u_h(t))v_h(t)\bigr)\bigr]\bigr|.
\end{split}
\end{equation*}
For all $t\in [t_0,t_0+\tau]$, let us define
\begin{equation*}
p_h(t)=M(u_h(t))v_h(t).
\end{equation*}
Denote by $\cald$ the set
\begin{equation*}
\begin{split}
\cald&=\bigl\{ (t,u_1, u_2, v, h): 0\le t \le T, \quad \vabs{u_1-u_0}\le
C_3\tau ,\\ &\quad \vabs{u_2-u_0}\le
C_3\tau,
\quad \vabs{v}\le C_3, \quad 0\le h \le h_1\}.
\end{split}
\end{equation*}
Let $L$ be
the Lipschitz constant of $(u_1,u_2)\mapsto F(t,u_1,u_2,v,h)$
restricted to $\cald$
and let $\omega_F$ be the modulus of continuity
of $F$ on $\cald$
With these notations, we can see that
\begin{equation}
\begin{split}
&\vabs{F^m- M(u_h(t))^{-1} f\bigl(t,u_h(t), M(u_h(t))\dot
u_h(t)\bigr)}\\
&\quad\le L\bigl(\vabs{U^m-u_h(t)}+\vabs{U^{m-1}-u_h(t)}\bigr)
\quad + 2\omega_F(h).\label{eq:84}
\end{split}
\end{equation}
Since $M$ is of class $C^3$ in $\Er^d$, $(u_h)_h$
converges strongly in $C^0([t_0,t_0+\tau])$ and $(v_h)_h$
converges strongly to $\dot u$ in $L^1(\Er)$ and almost everywhere on
$(t_0,t_0+\tau)$, the sequence $(p_h)_h$ also converges strongly
in $L^1(\Er)$ 
and almost everywhere on $(t_0,t_0+\tau)$ to $p=M(u)\dot u$. 
We see that  $M(u_h)^{-1}f(\cdot, u_h,p_h)$ tends to 
$M(u)^{-1}f(\cdot,u,p)$ strongly in $L^1(t_0,t_0+\tau)$ and almost
everywhere on $(t_0,t_0+\tau)$.
We summarize relations~\eqref{eq:85} and~\eqref{eq:84} together with
the above convergence result, and we find that
\begin{equation*}
\begin{split}
&\vabs{\langle F_h,\psi\rangle - \int_{t_0}^{t_0+\tau} \psi^T
M(u)^{-1}f(\cdot, 
u, M(u)\dot u)\, dt}\\
& \le
\int_{t_0}^{t_0+\tau} \vabs{M(u_h)^{-1}f\bigl(\cdot,
u_h,M(u_h)v_h\bigr)-M(u)^{-1}f\bigl(\cdot, u, M(u)\dot
u\bigr)}\vabs{\psi}\, dt\\ 
&\quad+ C_{8}\omega_\psi(h)\tau + (3LC_2 h +
2\omega_F(h))\int_0^T 
\vabs{\psi}\, 
dt,
\end{split}
\end{equation*}
which concludes the proof.
\end{proof} 

Let us prove now that the measure $\mu$ has the required variational
properties: 

\begin{lemma}\label{thr:1}
The measure $\mu$ satisfies
properties \eqref{eq:3}, \eqref{eq:4} and \eqref{eq:5}.
\end{lemma}

\begin{proof}Define
\begin{equation*}
\mu_h=M(u_h)\bigl(\ddot u_h -F_h\bigr);
\end{equation*}
$\mu_h$ is a sum of Dirac measures on $(t_0,t_0+\tau)$;
more precisely 
\begin{equation*}
\begin{split}
\mu_h&=\sum_{m=1}^nM(U^m)\bigl(V^m - V^{m-1}-hF^m\bigr)\delta(t-mh)\\
&\quad - M\bigl(U^{n+1}\bigr)V^n \delta(t-(n+1)h).
\end{split}
\end{equation*}
With all the previous results, we know that $\mu_h$ converges to
$\mu=M(u)\ddot u -f(\cdot, u,p)$ weakly in
$M^1\bigl((t_0,t_0+\tau)\bigr)$.  
Let us prove property~\eqref{eq:3}.
Assume that $\tau_0$ is a point of $(t_0,t_0+\tau)$ such
that $u(\tau_0)$ belongs 
to the interior of $K$. Then, by continuity of $u$, there exist
$\ve>0$ and $\rho>0$ such that
\begin{equation*}
\inf\{\vabs{u(t)-x}: \abs{t-\tau_0}\le \ve, x\in \partial K\}\ge 3\rho.
\end{equation*}
Since the sequence $\bigl(u_h\bigr)_h$ converges uniformly to $u$ as
$h$ tends to 
$0$, we can decrease $h_1$ so that
\begin{equation*}
\inf\{\vabs{u_h(t)-x}: \abs{t-\tau_0}\le \ve, 0<h\le h_1, x\in \partial
K\}\ge 2\rho.
\end{equation*}
Relation~\eqref{Wn} implies the identity
\begin{equation}
W^m=U^m +\frac{1-e}{1+e}hV^{m-1}+ \frac{h^2}{1+e} F^m.\label{eq:94}
\end{equation}
Relations~\eqref{eq:94} and~\eqref{eq:75} imply that
\begin{equation}
\vabs{W^m - U^m}\le hC_3\frac{1-e}{1+e}+ \frac{h^2C_{8}}{1+e}.\label{eq:179}
\end{equation}
Possibly decreasing $h_1$, we have thus
\begin{equation*}
\inf\{\vabs{W^m-x}: \abs{t_m-\tau_0}\le \ve, 0<h\le h_1, x\in \partial
K\}\ge \rho. 
\end{equation*}
This proves that the support of $\mu_h$ does not intersect the open
set $(\tau_0-\ve, \tau_0+\ve)$, and therefore, relation~\eqref{eq:3} holds.
Assume now that $u_1=u(t_1)$ belongs to $\partial K$, and let
$B(u_1,R_1)$ be a ball having the properties 
of theorem~\ref{thr:8}; assume that the image of  $(\tau_1, \tau_2)$
by $u_h$ and $w_h$ is included in this ball for all small enough
$h$.
We rewrite conditions~\eqref{eq:4} and~\eqref{eq:5} as follows: for
all continuous function $\psi$ with compact support included in $(t_0,
t_0+\tau)$ and taking its values in $\Er^d$ the following
implication holds:
\begin{equation}
\forall t\in (t_0,t_0+\tau), \quad d\phi(u(t))\psi(t)\ge
0\Longrightarrow 
\langle\mu, \psi\rangle\ge 0.\label{eq:86}
\end{equation}
In particular, if $d\phi(u(t))\psi(t)$ vanishes for all $t\in
(t_0,t_0+\tau)$, then $\langle \mu, \psi\rangle$ also vanishes.

The reader will check the equivalence of~\eqref{eq:4} and~\eqref{eq:5}
with~\eqref{eq:86}.
We infer from relation~\eqref{eq:169} that
\begin{equation*}
\vabs{Y(U^m)}\le \Lambda hC_3;
\end{equation*}
the above relation together with~\eqref{eq:179} imply that there
exists a constant $C_{10}$ such that
\begin{equation*}
\vabs{Y(W^m)}\le h C_{10}.
\end{equation*}
Since~\eqref{eq:86} is local, it is enough to check it in the
neighborhood of any $t_1\in (t_0, t_0+\tau)$. 
Let
\begin{equation*}
P=\lceil \tau_1/h\rceil, \quad Q=\lfloor \tau_2/h\rfloor,
\end{equation*}
and 
\begin{equation*}
\calp =\bigl\{m\in \{P,\dots,Q\}: W^m\notin K\bigr\}, \quad
\calp'=\{P,\dots,Q\} \setminus \calp.
\end{equation*}
We observe that if $m$ belongs to $\calp'$, then
\begin{equation*}
V^m -V^{m-1} - hF^m=0.
\end{equation*}
Therefore, we have the identity:
\begin{equation*}
\begin{split}
&\sum_{m=P}^Q \langle V^m - V^{m-1} - h F^m, \psi(t_m)\rangle_{U^m}\\
&\quad = \sum_{m\in \calp} \langle V^m - V^{m-1} - h F^m,
\psi(t_m)\rangle_{U^m}.
\end{split}
\end{equation*}
We recall relation~\eqref{diffVn}. Relation~\eqref{eq:180} implies that 
\begin{equation*}
\Phi(Z^m)-\Phi(W^m)=\begin{pmatrix}
0\\Y(W^m)^-
\end{pmatrix},
\end{equation*}
and therefore
\begin{equation*}
\vabs{Z^m - W^m -D\Psi(W^m)\begin{pmatrix}
0\\Y(W^m)^-
\end{pmatrix}} \le C_4 \vabs{Z^m -W^m}^2.
\end{equation*}
On the other hand, the definition of $\Psi$ is such that the $d$-th
column of $D\Psi(Z^m)$ is equal to $N(Z^m)$; therefore
\begin{equation*}
\begin{split}
&\vabs{D\Psi(W^m)\begin{pmatrix}
0\\Y(W^m)^-
\end{pmatrix} - N(Z^m) Y(W^m)^-}\\
&\quad\le 2C_4\vabs{Z^m - W^m} Y(W^m)^-.
\end{split}
\end{equation*}
We infer from the above estimates that
\begin{equation*}
\begin{split}
&\vabs{Z^m -W^m - Y(W^m)^-N(Z^m)}\\
&\quad\le C_4\bigl(2Y(W^m)^- +\vabs{Z^m -W^m}\bigr) \vabs{Z^m -W^m}\\
&\quad\le \frac{C_4 h}{1+e}\vabs{V^m - V^{m-1} -hF^m} \bigl(2+\Lambda\bigr)
h C_{10},
\end{split}
\end{equation*}
and thus, there exists $C_{11}$ such that for all $m\in \calp$:
\begin{equation*}
\vabs{Z^m - W^m -Y(W^m)^-N(Z^m)} \le h^2 C_{11} \vabs{V^m - V^{m-1} -hF^m}.
\end{equation*}
We can see now that
\begin{equation*}
\begin{split}
&\sum_{m\in\calp} \langle V^m - V^{m-1} - h F^m,
\psi(t_m)\rangle_{U^m} \\
&= \frac{1+e}{h} \sum_{m\in\calp} \langle Z^m
-W^m,\psi(t_m)\rangle_{U^m}\\
&\ge\frac{1+e}{h} \sum_{m\in\calp} Y(W^m)^- \langle N(Z^m),
\psi(t_m)\rangle_{U^m} \\
&\quad- C_{11} h(1+e)\max_{P\le m\le Q}
\bigl(\vnorm{M(U^m)}\vabs{\psi(t_m)} \sum_{m\in\calp} \vabs{V^m - 
V^{m-1} -h F^m},
\end{split}
\end{equation*}
which implies by a strightforward passage to the limit that $\langle
u, \psi\rangle$ is non negative. This concludes the proof of the lemma.
\end{proof}

\section{Transmission of energy during impact\label{sec:Transmission-energy}}

The basic assumption is still the one made at the beginning of
Section~\ref{sec:Estim-accel}. 

Let $\tab\in (0,\tau)$  be such that $u(\tab)$ belongs to $\partial
K$. Write $\tb=t_0+\tab$. We decompose
$p(\tb\pm0)$ into a normal component $p_N(\tb\pm0)$ 
belonging to $\Er d\phi(u(\tb))$ and a tangential part $p_T(\tb\pm 0)$
belonging to the orthogonal of $d\phi(u(\tb))$ in the cotangent metric
at $u(\tb)$.

In this section, we shall prove that
\begin{equation}
p_T(\tb+0)=p_T(\tb-0)\text{ and }
p_N(\tb+0)=-ep_N(\tb-0),\label{eq:95} 
\end{equation}
where $e$ is the restitution coefficient of the problem.

The conservation of the tangential component of the impulsion is
proved in next lemma:

\begin{lemma}\label{thr:3}
Assume that $\tab\in (0,\tau)$ is such that $u(\tab)$ belongs to
$\partial K$. Then 
\begin{equation*}
p_T(\tb+0)=p_T(\tb-0).
\end{equation*}
\end{lemma}

\begin{proof}
Thanks to lemma~\ref{thr:1}, we know that
\begin{equation}
M(u)\ddot u =\mu + f(\cdot, u, p),\label{eq:90}
\end{equation}
and that there exists a nonnegative measure $\lambda$ such that
\begin{equation}
\mu = \lambda d\phi(u).\label{eq:91}
\end{equation}
We take the measure of the set $\{\tb\}$ by the two sides
of~\eqref{eq:90}, and we find that
\begin{equation*}
M(u(\tb))\bigl(\dot u(\tb+0)-\dot
u(\tb-0)\bigr) = \mu(\{
\tb\}),
\end{equation*}
which implies immediately that $p(\tb+0)-p(\tb-0)$ is parallel to
$d\phi(u(\tb))$ and proves the lemma.
\end{proof}

Let $\overline{u}=u(\tb)$ and let
$B(\overline{u},r_1)$  and $B(\overline{u},r_1)$ have the properties
of theorem~\ref{thr:8}. There exists an 
interval $[\tau_{-5},\tau_{2}]$ containing $\tab$ in its interior such
the for all small enough $h$, $u_h\bigl([t_0+\tau_{-5},t_0+\tau_{2}]\bigr)$
is included in $B(u_1,r_1)$. 

The apparently strange notations $\tau_{-5}$ and $\tau_2$ have been
chosen in view of the upcoming construction of lemmas~\ref{thr:4}
and~\ref{thr:5}, where we will consider 
relative times
\begin{equation*}
\tau_{-5}< \dots<\tau_{-1}<\tab<\tau_1<\tau_2.
\end{equation*}

Define
\begin{equation*}
P=\lceil \tau_{-5}/h\rceil +1,\quad Q=\lfloor\tau_2/h\rfloor -1,
\end{equation*}
and let $x_h$ be obtained from the $X^m$ by affine interpolation, for
$P\le m\le Q$. We infer from estimates~\eqref{eq:168}
and~\eqref{eq:181} the estimates
\begin{equation*}
\begin{split}
&\max_{P\le m \le Q} \vabs{\frac{X^{m+1}-X^m}{h}}\le \Lambda
C_3,\\
&\sum_{m=P}^Q \vabs{\frac{X^{m+1}-X^m}{h} - \frac{X^m - X^{m-1}}{h}}
\le \Lambda C_7.
\end{split}
\end{equation*}
Therefore, we have the following convergences
\begin{equation*}
\begin{split}
x_h &\to x \text{ strongly in $C^0\bigl([t_0+\tau_{-5},t_0+\tau_2]\bigr)$;}\\
\dot x_h &\to \dot x \text{ except on a countable set and weakly
$*$}\\
&\qquad \qquad\qquad\text{in $ L^\infty\bigl([t_0+\tau_{-5},t_0+\tau_2]\bigr)$;}\\
\ddot x_h &\to \ddot x \text{ weakly in $M^1\bigl([t_0+\tau_{-5},t_0+\tau_2]\bigr)$.}
\end{split}
\end{equation*}
Write for all $h\le h_1$
\begin{equation*}
x_h= \begin{pmatrix}
s_h\\ y_h
\end{pmatrix}, \quad x= \begin{pmatrix}
s\\y
\end{pmatrix},
\end{equation*}
where the $s_h$'s and $s$ take their values in $\Er^{d-1}$ and the
$y_h$'s and $y$ are real valued functions.
We do not have $x_h=\Phi(u_h)$, because $x_h$ is a linear
interpolation  of the sequence $X^m=\Phi(U^m)$, and $\Phi(u_h)$ is the
image of the linear interpolation of the sequence $U^m$. However, we
can estimate the difference $x_h -\Phi(u_h)$.

\begin{lemma}\label{thr:11} For all $t \in [t_0+\tau_{-5},t_0+\tau_2]$,
belonging to $[t_m, t_{m+1}]$, we have:
\begin{equation*}
x_h (t) - \Phi\bigl(u_h(t) \bigr) \le 2C_4 C_3^2 h \min (t-t_m,
t_{m+1}-t).
\end{equation*}
\end{lemma}

\begin{proof}We observe that
\begin{equation*}
x_h(t_m)=X^m,
\end{equation*}
and that
\begin{equation*}
\begin{split}
&\vabs{\frac{d}{dt}
\bigl[x_h(t)-\Phi\bigl(u_h(t)\bigr)\bigr]\bigm|_{t=t_m+0}}\\
&\quad =\vabs{\frac{\Phi(U^{m+1})-\Phi(U^m)-hD\Phi(U^m)V^m}{h}}\\
&\quad \le hC_3^2C_{4}. 
\end{split}
\end{equation*}
Moreover, for all $t\in [t_m, t_{m+1})$
\begin{equation*}
\begin{split}
&\vabs{\frac{d^2}{dt^2}\biggl[x_h(t)-\Phi\bigl(u_h(t)\bigr)\biggr]
}\\
&\quad=\vabs{D^2\Phi\bigl(U^m +(t-t_m)V^m\bigr)V^m\otimes V^m}\le 2 C_3^2
C_{4}.
\end{split}
\end{equation*}
Therefore, a straightforward integration yields
\begin{equation*}
\vabs{x_h(t)-\Phi\bigl(u_h(t)\bigr)}\le C_3^2
C_{4}\bigl(h(t-t_m)+(t-t_m)^2\bigr), 
\end{equation*}
which implies
\begin{equation*}
\vabs{x_h(t)-\Phi\bigl(u_h(t)\bigr)}\le 2 C_3^2 C_{4} h(t-t_m).
\end{equation*}
We can write the analogous estimate on the interval $[t,t_{m+1})$,
which concludes the proof.
  \end{proof}

As a consequence of lemma~\ref{thr:11} we obtain:
\begin{equation*}
\forall t \in [t_0 +\tau_{-5}, t_0 +  \tau_2], \quad
x (t) = \Phi \bigl(u(t) \bigr), 
\end{equation*}
and
\begin{equation*}
\forall t \in (t_0 +\tau_{-5}, t_0 +  \tau_2), \quad
\dot x (t\pm 0) = D \Phi \bigl(u(t) \bigr) \dot u (t \pm 0). 
\end{equation*}
In virtue of relation~\eqref{eq:93},
\begin{equation*}
\dot u(\tb\pm 0)=\begin{pmatrix}
\dot s(\tb\pm 0)\\0
\end{pmatrix} + \dot y(\tb\pm 0) N(0).
\end{equation*}
We can rewrite this relation in terms of $p_N$ and $p_T$:
\begin{equation*}
p_T(\tb\pm 0)= M(0)\begin{pmatrix}
\dot s(\tb\pm 0)\\0
\end{pmatrix},\quad p_N(\tb\pm 0)= \dot y(\tb\pm 0) M(0)N(0)
\end{equation*}
Lemma~\ref{thr:3} implies $ \dot s(\tb + 0) =  \dot s(\tb -0)$.
In order to achieve the proof of relation~\eqref{eq:95}, we will prove
the scalar relation
\begin{equation}
\dot y(\tb +0)= -e \dot y(\tb -0).\label{4.35}
\end{equation}
We will do this by performing a precise analysis of the transmission
of energy on the scheme~\eqref{eq:106}.
The measure $\ddot y_h$ is a sum of Dirac measures on $(
t_0 +\tau_{-5}, t_0 +  \tau_2)$. We define two measures $\omega_h$ and
$\lambda_h$ on $(t_0 +\tau_{-5}, t_0 +  \tau_2)$ by  
\begin{equation*}
\omega_h =  \sum_{ m=P}^Q
\frac{\bigl(-2y^m + (1-e) y^{m-1} \bigr)^+}{ h} \delta
(t-mh),
\end{equation*}
and
\begin{equation*}
\lambda_h (t) = \sum_{m=P}^Q h
\lambda^m \delta (t-mh).
\end{equation*}
We have
\begin{equation*}
\ddot y_h=\omega_h +\lambda_h,
\end{equation*}
and it is obvious that
$\omega_h$ is a non-negative measure. 

Since the real numbers $\lambda^m$  are bounded independently of $h$
and $n$, the measure by  $\vabs{\lambda_h}$ of any subinterval $[a,b]$
of $(t_0+\tau_{-5}, t_0+\tau_2)$ is bounded by $C(b-a+h)$, and it is clear
therefore that there exists
a function $\lambda\in L^\infty(t_0+\tau_{-5}, t_0+\tau_2)$ and a subsequence
$\lambda_h$ converging to $\lambda$ in the weak topology of
$M^1\bigl((t_0+\tau_{-5}, t_0+\tau_2)\bigr)$. 

The measure $\omega_h$ converges in the weak topology of
$M^1\bigl((t_0+\tau_{-5}, t_0+\tau_2)\bigr)$ to a non-negative measure
$\omega$, and in the limit
\begin{equation}
\ddot y = \omega +\lambda,\label{eq:112}
\end{equation}
while
\begin{equation}
\abs[L^\infty]{\lambda}\le C_{9}.\label{eq:110}
\end{equation}
Since $y$ is non-negative on $(t_0+\tau_{-5}, t_0+\tau_2)$ and $y(\tau_0)$
vanishes, we must have
\begin{equation*}
\dot y(\tb+0)\ge 0, \quad \dot y(\tb-0)\le 0.
\end{equation*}
On  the other hand, $\dot y(\tb+0)  - \dot y(\tb-0)$ is equal to
$\omega(\{\tb\})$; if $\omega(\{\tb\})$ vanishes, we have
\begin{equation*}
\dot y(\tb+0)= \dot y(\tb-0)=0,
\end{equation*}
and the identity
\begin{equation*}
\dot y(\tb+0)=-e \dot y(\tb-0)
\end{equation*}
holds. Therefore, the only interesting case is when
\begin{equation}
\omega(\{\tb\})>0.\label{eq:133}
\end{equation}

The following two lemmas enable us to prove in two steps that the
velocity is reversed according to the law described
by~\eqref{eq:6}. Lemma~\ref{thr:4} shows that if $\omega$ has a Dirac
mass at $\tb$, then the left velocity at $\tb$ is outgoing;
Lemma~\ref{thr:5} shows indeed that~\eqref{eq:6} holds.

\begin{lemma}\label{thr:4} If $\omega(\{\tb\})$ is strictly positive, then $\dot
y(\tb-0)$ is strictly negative.
\end{lemma}

\begin{proof}The idea of the proof is to find two succesive times
$t_{m-1}\le t_m <\tb$ for which we can write down an estimate on
the discrete velocities, and then to use lemma~\ref{monodim} to
perform a discrete integration and to obtain a contradiction. We must
deal with the fact that $\dot y_h$ does not converge uniformly to
$\dot y$. 

Without loss of generality, we may assume that $\dot y$
is continuous on the right and that for all $h\le h_1$, $\dot y_h$ is
also continuous from the right. According to Helly's theorem, there
exists a countable set $D$ such that 
\begin{equation*}
\dot y_h(t) \to \dot y(t),\quad\forall t\text{ such that } t-\tb\in (\tau_{-5},
\tau_2)\setminus D.
\end{equation*}
Assume that $\dot y(\tb)$ vanishes; therefore, $\dot y(\tb+0)$
is strictly positive.
Choose $\alpha=\dot y(\tb+0)/4$, and let $\tau_{-4}$ and
$\tau_{1}$ be such that
\begin{align}
&\tau_{-5}\le \tau_{-4} < \tab<\tau_1 \le\tau_2\notag\\
&6C_{9} \bigl( \tau_1 -\tau_{-4}\bigr) \le \alpha,
\label{eq:115}\\
\intertext{and}
&\omega\bigl([t_0+\tau_{-4},\tb)\bigr)\le \alpha, \quad
\omega\bigl((\tb,t_0+\tau_{1}]\bigr)\le \alpha.
\end{align}
An integration of \eqref{eq:112} on appropriate intervals yields
\begin{align}
\forall t\in (t_0+\tau_{-4},\tb), \quad &\vabs{\dot y(t\pm 0)} \le
\alpha + C_{9} (\tb-t),\label{eq:113} \\
\forall t\in (\tb, t_0+\tau_1),\quad  &\dot y(t\pm 0) \ge
\omega\bigl(\{\tb\}\bigr) -\alpha - C_{9} (t-\tb).\label{eq:114}
\end{align}
Choose $\tau_{-3}\in (\tau_{-4}, \tab)\setminus D$ and $\tau_{-1}\in
(\tau_{-3},\tab)\setminus D$; since $\omega_h$ is a
nonnegative measure, we have the 
following inequality for all $\tau'\in (\tau_{-3},\tau_{-1})$ and all
$\tau''\in (\tau',\tau_{-1})$:
\begin{equation*}
\begin{split}
&\vabs{\dot y_h(t_0+\tau')-\dot y_h(t_0+\tau'')} \le \omega_h((t_0+\tau',t_0+\tau'']) +
C_{9}(\tau''-\tau' + h)\\
&\quad\le \omega_h([t_0+\tau_{-3},t_0+ \tau_{-1}]) +C_{9}(\tau''-\tau' + h).
\end{split}
\end{equation*}
We integrate $\omega_h - \omega$ on the interval
$\bigl[t_0+\tau_{-3},t_0+\tau_{-1}\bigr]$; since the measures $\omega$ and
$\omega_h$ do not charge $t_0+\tau_{-3}$ and
$t_0+\tau_{-1}$, we find that
\begin{equation*}
\begin{split}
&\omega_h\bigl(\bigl[t_0+\tau_{-3},t_0+ \tau_{-1}\bigr]\bigr)
- \omega\bigl(\bigl[t_0+\tau_{-3},t_0+ \tau_{-1}\bigr]\bigr)\\
&\quad= \dot y_h\bigl(t_0+\tau_{-1}\bigr) - \dot y_h\bigl(t_0+\tau_{-3}\bigr)
- \dot y\bigl(t_0+\tau_{-1}\bigr) + \dot y\bigl(t_0+\tau_{-3}\bigr)\\
& \quad+ \lambda
\bigl(\bigl[t_0+\tau_{-3}, t_0+\tau_{-1}\bigr]\bigr) -
\lambda_h\bigl(\bigl[t_0+\tau_{-3}, t_0+\tau_{-1}\bigr]\bigr),
\end{split}
\end{equation*}
and therefore
\begin{equation*}
\begin{split}
&\omega_h\bigl(\bigl[t_0+\tau_{-3}, t_0+\tau_{-1}\bigr]\bigr) \\
&\quad\le
\omega\bigl(\bigl[t_0+\tau_{-3}, t_0+\tau_{-1}\bigr]\bigr)  +\vabs{\dot
y_h(t_0+\tau_{-1})-\dot y(t_0+\tau_{-1})} \\
&\quad +\vabs{\dot
y_h(t_0+\tau_{-3})-\dot y(t_0+\tau_{-3})} + C_{9}\bigl(2\bigl(\tau_{-1} -
\tau_{-3}\bigr) + h\bigr).
\end{split}
\end{equation*}
Choose now $\tau_{-2}\in\bigl(\tau_{-3},\tau_{-1}\bigr) \setminus D$;
then, for 
$h$ small enough, $t_m=h\lfloor \tau_2/h\rfloor$ and  $t_{m-1}=t_m-h$
belong to the interval $(\tau_{-3}, \tau_{-1})$, and therefore, 
\begin{equation}
\vabs{\dot y_h(t_m)-\dot y_h(t_{m-1})}\le \alpha +
C_{9}\bigl(2\bigl(\tau_{-1}-\tau_{-3}\bigr) + 3h\bigr) +\ve_h,\label{eq:101}
\end{equation}
where $\ve_h$ tends to $0$ as $h$ tends to $0$.
On the other hand, $\dot y_h(t_0+\tau_{-2})$ tends to $\dot y(t_0+\tau_{-2})$
and therefore, thanks to relation~\eqref{eq:113}, 
there exists a family $\ve'_h$ such that
\begin{equation*}
\vabs{\dot y_h(t_0+\tau_{-2})}=\vabs{\dot y_h(t_m)} \le \alpha +
C_{9}\bigl(\tab-\tau_{-2}\bigr) + \ve'_h,
\end{equation*}
which is equivalent to
\begin{equation}
\vabs{\eta^m} \le \alpha +
C_{9}\bigl(\tab-\tau_{-2}\bigr) + \ve'_h;\label{eq:105}
\end{equation}
we infer from~\eqref{eq:101} and~\eqref{eq:105} that
\begin{equation*}
\vabs{\eta^{m-1}} \le 2\alpha +
C_{9}\bigl(2\bigl(\tau_{-1}-\tau_{-3}\bigr) 
+\tab-\tau_{-2}+ 3h\bigr) +\ve_h 
+\ve'_h.
\end{equation*}
Thus, for all $n \ge m$ we infer from Lemma~\ref{monodim} that
\begin{equation*}
\begin{split}
\vabs{\eta^m} &\le 2\alpha +
C_{9}\bigl(2\bigl(\tau_{-1}-\tau_{-3}\bigr) + 3h\\
&\quad  +\tab-\tau_{-2} +
2\bigl(t_m -t_m\bigr)\bigr) +\ve_h 
+\ve'_h .
\end{split}
\end{equation*}
Therefore, in the limit, for all $t\ge t_0+ \tau_{-2}$
\begin{equation*}
\vabs{\dot y(t)}\le 2\alpha + C_{9}\bigl(2\bigl(\tau_{-1}-\tau_{-3}\bigr)
+\tab-\tau_{-2} + 2\bigl(t-\tau_{-2}\bigr)\bigr).
\end{equation*}
and for all $t\in [t_0+\tau_{-2}, t_0+\tau_1]$
\begin{equation}
\vabs{\dot y(t)}\le 2\alpha + C_{9}\bigl(2\bigl(\tau_{-1}-\tau_{-3}\bigr)
+\tab-\tau_{-2} + 2\bigl(\tau_1-\tau_{-2}\bigr)\bigr).\label{eq:108}
\end{equation}
On the other hand, relation~\eqref{eq:114} implies that for all $t\in
(\tb, t_0+\tau_1)$,
\begin{equation}
\vabs{\dot y(t)} \ge 3\alpha - C_{9}
\bigl(\tau_{1}-\tab\bigr).\label{eq:109}
\end{equation}
Under assumption~\eqref{eq:115}, relation~\eqref{eq:109} contradicts relation~\eqref{eq:108}.
\end{proof}

We can conclude now the local study of the reflexion of the velocity
by the following lemma:

\begin{lemma}\label{thr:5}If $\omega\bigl(\{\tb\}\bigr)$ is
strictly positive, then 
\begin{equation}
\dot y(\tb)=-e\dot y(\tb-0).
\end{equation}
\end{lemma}

\begin{proof} Since $\dot y(\tb-0)$ is strictly negative, there
exists a real number $\tau_{-3}$ such that $y(t)$ is strictly positive
on $[t_0+\tau_{-3},\tb)\subset [t_0+\tau_{-5},\tb)$. For all
$\tau_{-2}\in (\tau_{-3},\tab)$, there exists $\tau_{-1}\in
(\tau_{-2}, \tab)$ and
$h_1>0$ such that
\begin{equation}
\forall h\in (0, h_1], \quad \forall  t \in [t_0+\tau_{-2},t_0+ \tau_{-1}),
y_h(t)\ge \frac{y(t_0+\tau_{-2})}{2}.\label{eq:126}
\end{equation}
We prove now that there exists a maximal integer 
\begin{equation*}
m\in
\{\lfloor\tau_{-3}/h\rfloor,\dots, \lfloor(\tau_0 + \ve)/h\rfloor\}
\end{equation*}
such
that
\begin{equation}
\forall l\in \bigl\{\lfloor \tau_{-3}/h\rfloor, \dots, m-1\bigr\}, \quad
2 y^l - (1-e)y^{l-1} \ge 0,\label{eq:127} 
\end{equation}
and denoting
\begin{equation}
\sigma_h = t_{m-1}-t_0,\label{eq:128}
\end{equation}
the time $\sigma_h$ satisfies
\begin{equation}
\lim_{h\to 0} \sigma_h =\tab.\label{eq:130}
\end{equation}
Let us first observe that for all small enough $h$ and all $t_l$
belonging to $[t_0+\tau_{-2},t_0+\tau_{-1}]$ we have
\begin{equation}
2y^l -(1-e)y^{l-1}\ge 0.\label{eq:129}
\end{equation}
Indeed,
\begin{equation*}
\begin{split}
2y^l -(1-e)y^{l-1} &= (1+e) y^l +(1-e)h \eta^{l-1} \\&\ge
\frac{1+e}{2} y(t_0+\tau_{-2}) - h (1-e) \Lambda C_{3},
\end{split}
\end{equation*}
and if 
$2 \Lambda C_{3} (1-e) h \le (1+e)y(t_0+\tau_{-2})$, we can see
that~\eqref{eq:129} holds. Therefore $m$ exists and
\begin{equation*}
\liminf \sigma_h \ge \tab.
\end{equation*}
On the other hand, if there existed $\tau_1>\tab$ such that for all
$t_m\in [t_0+\tau_{-3}, t_0+\tau_1]$ we had~\eqref{eq:129}, then
$\omega_h$ would vanish on $(t_0+\tau_{-3}, t_0+\tau_1)$, which
contradicts assumption~\eqref{eq:133}. Therefore, we have shown that
\begin{equation*}
\limsup \sigma_h \le \tab,
\end{equation*}
i.e.~\eqref{eq:130}.
We integrate discretely equation~\eqref{eq:154}, and we find that for
$t\in \bigl[t_0+\tau_{-3}, t_0+\sigma_h\bigr]$
\begin{equation}
\begin{split}
y_h(t)=&y_h(t_0+\sigma_h)-(t_0+\sigma_h-t) \dot y_h(t_0+\sigma_h)\\
&+\int_t^{t_0+\sigma_h} \lambda_h((s,t_0+\sigma_h])\, ds.
\end{split}
\label{eq:117}
\end{equation}
In the limit we have, 
\begin{equation}
y(t)=y(\tab)-(\tab-t)\lim_{h\downarrow 0} \dot y_h(t_0+\sigma_h+0) +
\int_{t}^{t_0+\tab} 
\int_{s}^{t_0+\tab} \lambda(r)\,dr\, ds.\label{eq:118}
\end{equation}
The comparison of~\eqref{eq:117} and~\eqref{eq:118} shows that
\begin{equation}
\lim_{h\downarrow 0}  \dot y_h(t_0+\sigma_h+0) = \lim_{h\downarrow 0}
 \eta^{m-1} =  \dot y(\tb-0).
\end{equation}

Our purpose now is to obtain very precise estimates on the behavior of
$y_h$ beyond $t_0+\sigma_h$.
Thanks to the maximality of of $m$, we have the relation
\begin{equation}
y^{m+1}=-e y^{m-1} + h^2\lambda^m;\label{eq:119}
\end{equation}
let us estimate $2y^{m+1} - (1-e)y^{m}$: we substitute the value of 
$y^{m+1}$ given by~\eqref{eq:119} into this expression, and we also
use~\eqref{eq:155} with  $m$ replaced by $m-1$; we find
\begin{equation*}
\begin{split}
&2y^{m+1} - (1-e)y^{m} \\
&\quad= -\bigl[2y^{m-1} -(1-e)y^{m-2}\bigr] -
(1-e)h^2\lambda^{m-1} + 2h^2\lambda^{m}.
\end{split}
\end{equation*}
We apply relation~\eqref{eq:106} for $n=m+1$ and we find that
\begin{equation*}
\begin{split}
&\eta^{m+1} + e \eta^{m-1} = h \bigl(\lambda^{m+1}-\lambda^m\bigr)\\
&\quad+ \bigl\{-\bigr[2y^{m-1} - (1-e)y^{m-2}\bigr]h^{-1} -
(1-e)h\lambda^{m-1} + 2h\lambda^m\bigr\}^+.
\end{split}
\end{equation*}
Therefore, we have
\begin{equation*}
\eta^{m+1} + e \eta^{m-1} \ge - 2 h C_{9}.
\end{equation*}
On the other hand, if $\xi= -\bigr[2y^{m-1} - (1-e)y^{m-2}\bigr]h^{-1} -
(1-e)h\lambda^{m-1} + 2h\lambda^m$ is lesser than or equal to $0$,
\begin{equation*}
\vabs{\eta^{m+1} + e \eta^{m-1}}\le 2 h C_{9};
\end{equation*}
if $\xi$ is positive, then the sign condition on $2y^{m-1} -
(1-e)y^{m-2}$ implies that
\begin{equation*}
\eta^{m+1} + e\eta^{m-1} \le h \bigl(\lambda^{m+1}+
\lambda^{m}\bigr) - (1-e) h \lambda^{m-1}.
\end{equation*}
Thus, we have shown that
\begin{equation}
\vabs{\eta^{m+1}+e\eta^{m-1}} \le 3 C_{9} h.\label{eq:124}
\end{equation}

If $e$ is strictly positive, then for all small enough $h$,
\begin{equation*}
\eta^{m+1}\ge e\vabs{\dot y(\tb-0)}/2.
\end{equation*}
Let us estimate now the expression $2y^{m+2} - (1-e)y^{m+1}$: we have
\begin{equation*}
2y^{m+2} -(1-e)y^{m+1}= -e\bigl[2 y^m -(1-e)y^{m-1}\bigr] + O(h^2).
\end{equation*}
If $2y^{m+2}-(1-e)y^m$ is non-negative, then
\begin{equation*}
y^{m+3}=2y^{m+2}-y^{m+1} + h^2\lambda^{m+2}.
\end{equation*}
We must estimate $2y^{m+3} -(1-e)y^{m+2}$:
\begin{equation*}
\begin{split}
&2y^{m+3} -(1-e)y^{m+2} - 2y^{m+2}+(1-e)y^{m+1}\\
&=
h\bigl(2\eta^{m+2} - (1-e)\eta^{m+1}\bigr)\\
&=h(1+e)\eta^{m+1} +2h^2\lambda^{m+2},
\end{split}
\end{equation*}
and therefore $2y^{m+3} -(1-e)y^{m+2}$ is non negative for all small
enough $h$; the repetition of the argument shows that there exists
$\theta>0$ such that for all small enough $h$ and all
$n\in\{m+2,\dots, m+\lfloor\theta/h\rfloor\}$, the
expression $2y^{m+1} -(1-e)y^m$ is non negative, and thus we have the
relations
\begin{equation*}
y^{m} = y^{m+1} + h(n-m-1)\eta^{m+1} + \sum_{j=m+2}^{m-1} \bigl(n-j)
h^2 \lambda^j.
\end{equation*}

On the other hand, if $2y^{m+2} -(1-e)y^m$ is negative, we must
have
\begin{equation*}
y^m=-\frac{(1-e)h\eta^{m-1}}{1+e} + O(h^2),
\end{equation*}
and therefore
\begin{equation*}
y^{m-1} =-\frac{2h\eta^{m-1}}{1+e} + O(h^2).
\end{equation*}
These relations and the assumption on the sign of $2y^{m+2}
-(1-e)y^m$ imply that
\begin{equation}
2y^{m+3} -(1-e)y^{m+2} =-\frac{\bigl(4e^2
+e(1-e)^2\bigr)h\eta^{m-1}}{1+e} + O(h^2),\label{eq:125}
\end{equation}
which is strictly positive for $h$ small enough. But now, we can see
that
\begin{equation*}
y^{m+3}-y^{m+2} = -e h\eta^{m-1} + O(h^2),
\end{equation*}
which is strictly positive for small enough $h$, and therefore
$2y^{m+4}-(1-e)y^{m+3}$ is strictly positive for $h$ small enough,
since
\begin{equation*}
2y^{m+4}-(1-e)y^{m+3} \ge -he(1+e)\eta^{m-1} + O(h^2);
\end{equation*}
the same argument as above shows now that there exists $\theta>0$ such
that for all $n\in\{m+3,\dots,m+\lfloor\theta/h\rfloor\}$, 
\begin{equation*}
y^{m} = y^{m+2} + h(n-m-2)\eta^{m+2} + \sum_{j=m+3}^{m-1} \bigl(n-j)
h^2 \lambda^j.
\end{equation*}
If we let $\sigma'_h =t_{m+1}-t_0$ in the first case and
$\sigma'_h=t_{m+2}-t_0$ in the second case, we have now for 
$\sigma'_h\le t-t_0\le \sigma'_h+\theta-h$
\begin{equation}
y_h(t) = y_h(t_0+\sigma'_h) + (t-\sigma'_h-t_0)\dot y_h(t_0+\sigma'_h) +
\int_{t_0+\sigma'_h}^t \lambda_h((s,t])\, ds\label{eq:122}
\end{equation}
and 
\begin{equation}
y_h(t_0+\sigma'_h)=O(h), \quad \dot y_h(t_0+\sigma'_h)=-e \eta^{m-1} +
O(h).\label{eq:123}
\end{equation}
Passing to the limit in~\eqref{eq:122}, we can see that
\begin{equation*}
\dot y(\tb+0)=-e\dot y(\tb-0).
\end{equation*}

If we assume now that $e$ vanishes, relation~\eqref{eq:124} implies
\begin{equation*}
\eta^{m+1}=0(h).
\end{equation*}
We observe that lemma~\ref{monodim} implies that for all $n$
\begin{equation*}
\vabs{\eta^m} \le \vabs{\eta^{m-1}} + 2 C_{9} h,
\end{equation*}
which implies
immediately that for $n\ge m+1$
\begin{equation*}
\vabs{\eta^{m}} \le \vabs{\eta^{m+1}} + 2hC_{9}(n-m-1),
\end{equation*}
which proves by a straightforward passage to the limit that
\begin{equation*}
\dot y(\tb+0)=0.
\end{equation*}
This completes the proof of the lemma.
\end{proof}

\section{Initial conditions}\label{sec:Initial-conditions}

In this section we prove that the solution that we have constructed
satisfies the initial conditions; we work under the hypotheses stated
at the beginning of section~\ref{sec:Estim-accel}.

\begin{lemma}
The function $u$ satisfies the initial conditions
\begin{equation*}
u(0)=u_0, \quad p(0+0)=p_0.
\end{equation*}
\end{lemma}

\begin{proof}
By uniform convergence of $u_h$ to $u$, it is clear that $u(0)$ is
equal to $u_0$. There remains to show that the initial condition on
the impulsion is satisfied.

Assume first that $u_0$ belongs to the interior of $K$; then there
exist $h_1>0$ and $\tau_1>0$ such that for all $h\in (0, h_1]$ and for
all $t\in [0,\tau_1]$
\begin{equation*}
\vabs{u_h(t)-u_0}\le \frac{1}{2} \inf\{\vabs{u_0-y}: y\notin K\}.
\end{equation*}
Then for all $t_m$ belonging to $(0,\tau_1]$,
$\bigl(2U^{m} - (1-e)U^{m-1} + h^2F^m\bigr)/(1+e)$ belongs to $K$
for $h$ small enough; we have indeed
\begin{equation*}
\begin{split}
&\vabs{\frac{2U^{m} - (1-e)U^{m-1} + h^2F^m}{1+e} - u_0}\\
&\quad \le \frac{1-e}{1+e}hC_{1} +  \frac{1}{2} \inf\{\vabs{u_0-y}:
y\notin K\} + 
\frac{h^2}{1+e}C_{8},
\end{split}
\end{equation*}
which is strictly inferior to $\inf\{\vabs{u_0-y}: y\notin K\}$ for
$h$ small enough. Thus the constraints are not saturated for $0\le t_m
\le \tau_1$ and the convergence is a classical result.

In the second case, $u_0$ belongs to $\partial K$; we have taken
admissible initial conditions, so that
\begin{equation*}
\langle p_0, d\phi(u_0)\rangle_{u_0}^*\ge 0.
\end{equation*}
We use the construction and notations of section~\ref{sec:Existence}: $\Phi$,
$\Psi$, $X^m$, $s^m$, $y^m$ and $\zeta^m$ have the same signification
as there.

Taylor's formula yields
\begin{equation*}
\zeta^0 =\frac{X^1 - X^0}{h} = D\Phi(u_0)\frac{U^1 - u_0}{h} + O(h),
\end{equation*}
and the definition~\eqref{u1} of $U^1$ gives
\begin{equation}
\zeta^0 = D\Phi(u_0)M(u_0)^{-1} p_0 + O(h).\label{eq:140}
\end{equation}
Write
\begin{equation*}
\begin{pmatrix}
\sigma_0\\\eta_0
\end{pmatrix}=D\Phi(u_0)M(u_0)^{-1} p_0.
\end{equation*}
Then the normal and tangential components of the impulsion are given
by
\begin{equation*}
p_{0T}= M(u_0)\begin{pmatrix}
\sigma_0\\0
\end{pmatrix} \text{ and } p_{0N}=\eta_0 M(u_0) N(u_0).
\end{equation*}
We wish to prove
\begin{equation*}
p(0+0)=p_0,
\end{equation*}
which is equivalent to
\begin{equation*}
\dot x(0+0)=\begin{pmatrix}
\sigma_0\\ \eta_0
\end{pmatrix}.
\end{equation*}
We recall relation~\eqref{eq:154}.
Relation~\eqref{eq:140} implies that
\begin{equation*}
\sigma^1 =\bigl(D\Phi(u_0)M(u_0)^{-1}p_0\bigr)' +O(h),
\end{equation*}
and together with~\eqref{eq:154}, we obtain in the limit
\begin{equation*}
\dot s(t)=\bigl(D\Phi(u_0)M(u_0)^{-1}p_0\bigr)' + O(t),
\end{equation*}
i. e.
\begin{equation*}
\dot s(0+0)= \sigma_0.
\end{equation*}
Let us show now that
\begin{equation*}
\dot y(0+0)=\eta_0, 
\end{equation*}
considering two cases: $\eta_0>0$ and $\eta_0=0$. When $\eta_0$
vanishes, we have
\begin{equation*}
y^1= y^0 + h\eta_0 + O(h^2)= O(h^2),
\end{equation*}
and
\begin{equation*}
\begin{split}
y^2&= -ey^0 + \bigl(2y^1 -(1-e)y^0\bigr)^+ + h^2 \lambda^1\\
&= 2\bigl(y^1\bigr)^+ + h^2\lambda^1=O(h^2).
\end{split}
\end{equation*}
Thus,
\begin{equation*}
\eta^0= O(h),\quad \eta^1=O(h),
\end{equation*}
and relation~\eqref{eq:131} implies
\begin{equation*}
\vabs{\eta^m}\le O(h)+2 C_{9} h(n-1);
\end{equation*}
therefore, a passage to the limit gives immediately
\begin{equation*}
\dot y(0+0)=0.
\end{equation*}
If, on the other hand, $\eta^0$ is strictly positive, then
\begin{equation*}
2 y^1 -(1-e)y^0= 2y^1 = 2h\eta^0 + O(h^2)
\end{equation*}
which is strictly positive if $h$ is small enough. Let $\{1,\dots,m\}$
be the maximal interval such that
\begin{equation*}
2 y^n -(1-e)y^{n-1} >0, \quad \text{ if } n\le m.
\end{equation*}
Then, for all $n\in \{1,\dots, m\}$,
\begin{equation*}
\eta^n - \eta^{n-1} = h\lambda^{n},
\end{equation*}
which implies by discrete integration that
\begin{equation*}
\eta^n \ge \eta_0 - hn C_{9},
\end{equation*}
as long as $n$ belongs to $\{1,\dots, m\}$.
Moreover, if we choose any $\tau_1<\eta_0/(2C_{9})$ and if $n$ is at
most equal to $\min\bigl(m,\lfloor \tau_1/h\rfloor\bigr)$, we can see that
\begin{equation*}
y^m = y_0 + h\bigl(\eta^0+\dots+\eta^{m-1}\bigr) \ge \frac{hn\eta_0}{2},
\end{equation*}
for all small enough values of $h$. 

In particular, for all $n\le
\min \bigl(m,\lfloor \tau_1/h\rfloor\bigr)$,
\begin{equation*}
2 y^m -(1-e) y^{m-1} \ge \frac{(1+e) hn\eta_0}{2} - (1-e)h \Lambda C_{3},
\end{equation*}
which proves that $m$ is at least equal to
$\lfloor\tau_1/h\rfloor$. Therefore, $\omega_h$ vanishes on the
interval $(0, \tau_1-h)$; in the limit, $\omega$ vanishes on
$(0,\tau_1)$ and therefore
\begin{equation*}
\dot y(0)=\eta_0,
\end{equation*}
which completes the proof of the lemma.
\end{proof}

\section{A priori estimates}\label{sec:apriori}

In this section we prove that solutions of the problem \eqref{eq:2},
\eqref{eq:3}, \eqref{eq:4},\eqref{eq:5},
\eqref{eq:7}, \eqref{eq:8}, \eqref{position} and
\eqref{eq:9} satisfy an a priori estimate on an interval with non
empty interior.

\begin{lemma}\label{thr:6}
Let $R$ be strictly larger than $\abs[u_0]{p_0}^*$. Then
there exists $\tau(R)>0$ such that for all solution $u$ of \eqref{eq:2},
\eqref{eq:3}, \eqref{eq:4},\eqref{eq:5},
\eqref{eq:7}, \eqref{eq:8}, \eqref{position} and
\eqref{eq:9} defined on $[t_0,t_0+\tau]$, the following estimates hold:
\begin{equation}
\forall t \in [t_0,t_0+\min(\tau, \tau(R))], \quad 
\abs{u(t)-u_0}\le R, \quad \abs[u(t)]{p(t)}^*\le R.\label{eq:134}
\end{equation}
\end{lemma}

\begin{proof} The measure $\lambda$ appearing in~\eqref{eq:4} can be
decomposed in the sum of an atomic part $\lambda_a$ and a diffuse part
$\lambda_d$. At each point of the support of $\lambda_a$ we have
\begin{equation}
\abs[u(t)]{p(t+0)}^* \le \abs[u(t)]{p(t-0)}^*\label{eq:10}
\end{equation}
thanks to relation~\eqref{eq:6}. On any interval $(t_1,t_2)$ which
does not intersect the support of $\lambda_a$, we multiply
relation~\eqref{eq:2} by $\dot u^T$ on the left, and we find that
\begin{equation}
\frac{d}{dt} \frac{1}{2}\dot u^T M(u)\dot u = \dot u^T f(\cdot,u,p)
+ \frac{1}{2} \dot u^T (DM(u)\dot u)\dot u.\label{eq:11}
\end{equation}
Define
\begin{equation*}
E(u, p)= \frac{1}{2}\langle p,p\rangle_u^*, \quad z=\abs[u]{p}^*.
\end{equation*}
It is convenient to recall that
\begin{equation*}
\abs[u]{p}^*=\abs{M(u)^{-1/2}p}=\abs{M(u)^{1/2}\dot u}.
\end{equation*}
Relations~\eqref{eq:10} and~\eqref{eq:11} imply that in the sense of
measures
\begin{equation}
z \dot z = \dot E \le \dot u^T f(\cdot, u,p)+ \frac{1}{2} \dot u^T
\bigl(DM(u)\dot u\bigr)\dot u.\label{eq:12}
\end{equation}
Our purpose now is to transform~\eqref{eq:12} into a differential
inequality. Let $\chi(u)$ be the norm of the bilinear mapping
\begin{equation*}
(v_1, v_2)\mapsto M(u)^{-1/2}(DM(u)M(u)^{-1/2}v_1)M(u)^{-1/2}v_2.
\end{equation*}
With this definition,
\begin{equation*}
\vabs{\dot u^T\bigl(DM(u)\dot u\bigr)\dot u}\le \chi(u)z^3.
\end{equation*}
We write now
\begin{equation*}
\begin{split}
\dot u^T f(t,u,p)&=\dot u^TM(u)^{1/2}M(u)^{-1/2}f(t,u,p)\\
&=\dot u^TM(u)^{1/2} \bigl[M(u)^{-1/2}f(t,u,p)-M(u_0)^{-1/2}f(t,u_0,
p_0) \\
&+ M(u_0)^{-1/2} f(t,u_0, p_0)\bigr].
\end{split}
\end{equation*}
Define
\begin{equation*}
g(t)=\abs{M(u_0)^{-1/2}f(t,u_0,p_0)},
\end{equation*}
and let $\omega(\tau,R)$ be the Lipschitz constant of $(u,p)\mapsto
M(u)^{-1/2} f(t,u,p)$ for $t\in [t_0,t_0+\tau]$ and $\max(\abs{u-u_0},
\abs[u]{p}^*)\le R$; more precisely
\begin{equation*}
\begin{split}
\omega(\tau, R)&=
\sup\biggl\{
\frac{
\abs{M(u_1)^{-1/2}
f(t,u_1,p_1)
-{M(u_2)}^{-1/2}
f(t,u_2,p_2)}
}
{\abs{u_1-u_2}+\abs{p_1-p_2}}
:\\
& t_0\le t \le
t_0+\tau, \max(\abs{u_1-u_0},\abs{u_2-u_0}, \abs[u_1]{p_1}^*,
\abs[u_2]{p_2}^*)\le 
R, \\
&u_1\neq u_2\text{ or }  p_1\neq p_2\biggr\}.
\end{split}
\end{equation*}
By construction, $\omega$ is continuous and it is an increasing
function of $\tau$ and $R$.

Fix $R>\abs[u_0]{p_0}^*$.

If $t_0\le t \le t_0+\tau$ and if
$\max(\abs{u(t)-u_0},\abs[u(t)]{p(t)}^*)\le R$ on $[t_0, t_0+\tau]$,
we have the 
inequality
\begin{equation*}
\vabs{\dot u^T f(\cdot,u,p)} \le z\bigl(g+\omega(\tau,R)(\abs{u-u_0}
+ \abs{p-p_0})\bigr).
\end{equation*}
But we can estimate $u(t)-u_0$:
\begin{equation*}
\abs{u(t)-u_0}\le\int_{t_0}^t \abs{\dot u(s)}\, ds \le \int_{t_0}^t
\norm{M(u)^{-1/2}} z\, ds.
\end{equation*}
Therefore we have the estimate
\begin{equation*}
\begin{split}
\abs{\dot u^T f(\cdot,u,p)}&\le zg+z \omega(\tau,R)\biggl(\int_{t_0}^t
\norm{M(u)^{-1/2}}z \, ds\\
&+\norm{M(u)^{1/2}}z +\abs{p_0}\biggr),
\end{split}
\end{equation*}
and we conclude that $z$ satisfies the differential inequality
\begin{equation*}
\dot z \le g + \omega(\tau,R)\Bigl[\int_{t_0}^t \norm{M(u)^{-1/2}}z \, ds
+ \norm{M(u)^{1/2}}z + \abs{p_0}\Bigr] +\frac12 \chi(u) z^2.
\end{equation*}

Set
\begin{align}
\alpha&=\sup\bigl\{ \norm{M(u)^{1/2}}: \abs{u-u_0}\le R\bigr\},\label{eq:18}\\
\beta&=\sup\bigl\{ \norm{M(u)^{-1/2}}: \abs{u-u_0}\le R\bigr\},\label{eq:19}\\
\gamma&=2\sup\bigl\{\chi(u): \abs{u-u_0}\le R\bigr\}.\notag
\end{align}
While $t\le t_0+ \tau$ and $\max(\abs{u(t)-u_0}, \abs[u(t)]{p(t)}^*)\le R$,
$z$ satisfies the following differential inequality
\begin{equation}
\dot z \le g +\omega(\tau, R)\Bigl[\beta \int_{t_0}^t z\, ds +\alpha z +
\abs{p_0} \Bigr] + \gamma z^2.\label{eq:14}
\end{equation}

Let $\rho$ be any positive number; consider the integrodifferential equation
\begin{equation}
\dot y=g + \rho\Bigl(\beta \int_{t_0}^t y\, ds + \alpha y +\abs{p_0}\Bigr)
+ \gamma\abs{y}^2,\label{eq:15}
\end{equation}
with the initial condition
\begin{equation*}
y({t_0})=z({t_0}).
\end{equation*}
It has a unique maximal solution which blows up in finite time, as
soon as $\gamma$ is strictly positive and $\rho\abs{p_0}+ \sup\abs{g}$ is
strictly positive. 
Let $\theta(\rho)\in [{t_0},T]$ be the largest time for which
\begin{equation*}
\forall t \in [{t_0}, \theta(\rho)], \quad y(t)\le R,\quad \beta \int_{t_0}^t
y\, ds \le R.
\end{equation*}
As $\theta$ is a decreasing function of $\rho$, there exists a unique
$\tau(R)$ such that
\begin{equation*}
\theta(\omega(\tau(R),R))=\tau(R).
\end{equation*}
Choose now 
\begin{equation*}
\rho=\omega(\tau(R),R).
\end{equation*}
Then we can compare the solution $z$ of~\eqref{eq:14} and the solution
$y$ of~\eqref{eq:15}, and we find immediately that
\begin{equation}
\forall t\in \bigl[{t_0},t_0+\min\bigl(\tau,\tau(R)\bigr)\bigr],\quad
z(t)\le y(t).\label{eq:13} 
\end{equation}
This concludes the proof of the lemma.
\end{proof}

\section{Global results}\label{sec:Global}

We summarize the results obtained so far in the following Proposition:

\begin{proposition}\label{recap}Assume that there
exist strictly positive numbers $\tau$, $C_3$ and $h_1>0$, and a
sequence of solutions of the numerical scheme defined
by~\eqref{u0}, \eqref{u1}, \eqref{un} and~\eqref{fn}, which 
satisfies the estimate~\eqref{eq:168}. Then it is
possible to extract from the sequence $u_h$ defined by~\eqref{eq:79}
a subsequence which converges to a solution
of \eqref{eq:2},
\eqref{eq:3}, \eqref{eq:4},\eqref{eq:5},
\eqref{eq:7}, \eqref{eq:8}, \eqref{position} and \eqref{eq:9}. The
convergence holds in the following sense: $u_h$ converges uniformly to
$u_h$ on $[t_0, t_0+\tau]$; $\dot u_h$ converges to $\dot u$ in
$L^\infty(t_0,t_0+\tau)$ weaky star and almost everywhere on $[t_0,t_0+\tau]$,
and $\ddot u_h$ converges to $\ddot u$ in the weak topology of
measures.  Moreover, for all $\tau\in (t_0,t_0+\tau]$, we have the
following convergence:
\begin{equation}
\begin{split}
&\limsup_{h\downarrow 0} \sup\bigl\{\vabs[U^m]{V^m}: t_0\le t \le
t_0+\tau\bigr\} \\
&\quad\le \esssup\bigl\{ \vabs[u(t)]{\dot u(t)}: t_0\le t \le
t_0+\tau\bigr\}.\label{eq:136} 
\end{split}
\end{equation}
\end{proposition}

\begin{proof}
The only statement which deserves a proof is the last one; if it is
not true, there exists $\gamma >0$, a sequence of time steps still
denoted by $h$ and a sequence of integers $m(h)$ such that
\begin{equation}
\vabs[U^{m(h)}]{V^{m(h)}}^2\ge \esssup \bigl\{ \vabs[u(t)]{\dot
u(t)}^2: t_0\le t \le t_0+\tau\bigr\} +\gamma.\label{eq:43}
\end{equation}
Without loss of generality, we may assume that $h m(h)$ tends to
$\tau_2\in[0,\tau]$.

First, $\tau_2$ cannot be equal to $0$: we have learnt in
section~\ref{sec:Initial-conditions} that there exists a constant
$C_{12}$ and a time $\tau_1$ such that for all $h\le h_1$ and all
$m\le \tau_1/h$,
\begin{equation*}
\vabs{V^m -V^0}\le C_{12} mh.
\end{equation*}
In particular, this estimate implies that
\begin{equation*}
\vabs[U^{m(h)}]{V^{m(h)}} = \vabs[u_0]{v_0} + O(mh);
\end{equation*}
but $\vabs[u_0]{v_0}$ is at most equal to 
$\bigl\{ \vabs[u(t)]{\dot u(t)}: t_0\le t \le
t_0+\tau\bigr\}$, which contradicts~\eqref{eq:43}.
In the same fashion, we cannot have $u(t_0+\tau_2)\in \inter{K}$; if
it were the case, we could find an interval $[\tau_1, \tau_3]$
containing $\tau_2$ and $h_1>0$ such that for all $h\in (0, h_1]$,
 $u_h([\tau_1, \tau_3])$ is included in a ball of radius $r$ about
$u(t_0+\tau_2)$ included in the interior of $K$. But, in this case,
$\dot u_h$ converges uniformly to $\dot u$ in $C^0([\tau_1,\tau_3])$
and this contradicts again~\eqref{eq:43}.

Thus, we assume that $\tau_2$ is strictly positive and that
$u(t_0+\tau_2)$ belongs to $\partial K$. Choose a coordinate system
such that the origin is at $u(t_0+\tau_2)$; let $\Psi$ be the
diffeomorphism defined at~\eqref{eq:44}. In this case, $D\Psi(0)$ is
given by~\eqref{eq:93}. Define
\begin{equation*}
\beta^m=\bigl(\xi^m\bigr)^T D\Psi(0)^T M(0)D\Psi(0)\xi^m.
\end{equation*}
Let us compare $\beta^m$ to $\vabs[U^m]{V^m}^2$; it is convenient to
define
\begin{equation*}
\tilde V^m=D\Psi(0)\xi^m;
\end{equation*}
then
\begin{equation*}
\begin{split}
&\vabs[U^m]{V^m}^2-\beta^m =\bigl(V^m\bigr)^TM(U^m)V^m -\bigl(\tilde
V^m\bigr)^T M(0)\tilde V^m\\
&\quad=\bigl(V^m\bigr)^T\bigl(M(U^m)-M(0)\bigr)V^m +\bigl(V^m-\tilde
V^m\bigr) M(0)\bigl(V^m -\tilde V^m)\\
&\quad + 2\bigl(V^m -\tilde V^m\bigr)^T
M(0)V^m.
\end{split}
\end{equation*}
We observe that
\begin{equation*}
\vabs{U^m}\le \vnorm{u_h-u} + C_3\vabs{mh-\tau_1},\quad \vnorm{X^m}\le
\vnorm{U^m},
\end{equation*}
and that
\begin{equation*}
\vabs{V^m -\tilde V^m} \le C_4\vnorm{\xi^m}\bigl[2\Lambda \vnorm{X^m}
+\vnorm{X^m -X^{m-1}}\bigr].
\end{equation*}
These observations enable us to estimate the difference: there exists
a constant $C_{13}$ such that
\begin{equation*}
\begin{split}
\vabs{\vabs[U^m]{V^m}^2 - \beta^m}
\le C_{13}\bigl(h + \vnorm[{C^0([t_0,t_0+\tau])}]{u -u_h} 
+\vabs{mh -\tau_1}\bigr).
\end{split}
\end{equation*}
We infer from~\eqref{eq:51} that there exists a constant $C_{14}$ such that
\begin{equation*}
\beta^{m+1}\le \min\bigl(\beta^m,\beta^{m-1}\bigr) +C_{14}h.
\end{equation*}
We use now~\eqref{eq:43}: we can see that for all $m\le m(h)$,
\begin{equation*}
\beta^{m(h)} \le\max\bigl(\beta^m, \beta^{m-1}\bigr) + C_{14}\bigl(m(h)-m\bigr),
\end{equation*}
so that
\begin{equation*}
\begin{split}
&\max\bigl(\vabs[U^{m}]{V^m}, \vabs[U^{m-1}]{V^{m-1}}\bigr) \ge
\beta^{m(h)} - C_{14}\bigl(m(h)-m\bigr)h -\\
&\qquad C_{13}\bigl(h + \vnorm[{C^0([t_0,t_0+\tau])}]{u -u_h} 
+\vabs{mh -\tau_1}\bigr).
\end{split}
\end{equation*}
If $\tau_4<\tau_1$ is such that
\begin{equation*}
\tau_1 -\tau_4 \le \gamma/(4C_{14}),
\end{equation*}
and if
\begin{equation*}
C_{13}\bigl(h + \vnorm[{C^0([t_0,t_0+\tau])}]{u -u_h} 
+\vabs{mh -\tau_1}\bigr)\le \gamma/(4C_{14}),
\end{equation*}
we can see that for all small enough $h$ and all $m\in\{\lceil
\tau_4/h\rceil,\dots, m(h)\}$ the following estimate holds:
\begin{equation}
\begin{split}
&\max\bigl(\vabs[U^{m}]{V^m}, \vabs[U^{m-1}]{V^{m-1}}\bigr)\\
&\qquad\ge 
\esssup \bigl\{ \vabs[u(t)]{\dot
u(t)}^2: t_0\le t \le t_0+\tau\bigr\}
+\gamma/2.
\end{split}\label{eq:46}
\end{equation}
But the function $v_h$ defined by
\begin{equation*}
v_h(t)= \vabs[U^m]{V^m}^2 \text{  if } t \in [mh, (m+1)h)
\end{equation*}
converges almost everywhere on $[0,\tau]$ to $\vabs[u(t)]{\dot
u(t)}^2$; so does $\max(v_h(t-h),v_h(t))$. Therefore, in the
limit, relation~\eqref{eq:46} leads to
\begin{equation*}
\esssup_{t\in [t_0+\tau_{4},t_0+\tau_{2}]}v_h(t) \ge \esssup \bigl\{ \vabs[u(t)]{\dot
u(t)}^2: t_0\le t \le t_0+\tau\bigr\}+ \gamma/2,
\end{equation*}
which is a contradiction.
\end{proof}

A corollary can be inferred imediately from this Proposition and
Theorem~\ref{thr:12}:

\begin{corollary}\label{thr:7}
For all admissible initial conditions $u_0$ and $p_0$, there exists
$\tau>0$ and a solution of \eqref{eq:2},
\eqref{eq:3}, \eqref{eq:4},\eqref{eq:5},
\eqref{eq:7}, \eqref{eq:8}, \eqref{position} and \eqref{eq:9}
defined on $[t_0,t_0+\tau]$.
\end{corollary}

We have proved above the existence of a non-empty interval on which
the numerical scheme converges to a solution of \eqref{eq:2},
\eqref{eq:3}, \eqref{eq:4},\eqref{eq:5},
\eqref{eq:7}, \eqref{eq:8}, \eqref{position} and \eqref{eq:9}. 
On the other hand, lem\-ma~\ref{thr:6} gives a priori estimates on the
solution of such a problem.

We couple now the a priori estimates with the local
convergence result to obtain a global result:

\begin{theorem}Let $R$ be strictly larger than $\vabs[u_0]{p_0}^*$,
and let $\tau(R)$ be given as in lemma~\ref{thr:6}. Then, for all small
enough $h$, the solution $U^m$ of the numerical scheme~\eqref{u0},
\eqref{u1}, \eqref{un}, \eqref{fn} is defined on a discrete interval
$\{0,\dots,m(h)\}$, such that
\begin{equation*}
hm(h) \to \tau(R);
\end{equation*}
moreover, the approximation $u_h$ converges to a solution $u$ of the
continuous time equation, i.e.
\eqref{eq:2}, 
\eqref{eq:3}, \eqref{eq:4},\eqref{eq:5},
\eqref{eq:7}, \eqref{eq:8}, \eqref{position} and \eqref{eq:9}, which
is defined on $[t_0, t_0+\tau(R)]$.
\end{theorem}

\begin{proof} Let 
Let $\{0,\dots,m(h)\}$ be the discrete time interval for which the
numerical scheme ~\eqref{u0},
\eqref{u1}, \eqref{un}, \eqref{fn} has a solution; we know from
theorem~\ref{thr:12} that
\begin{equation*}
\liminf h m(h)=\tau>0.
\end{equation*}
Assume that
\begin{equation}
\tau <\tau(R).\label{eq:47}
\end{equation}
It is possible to extract from the sequence $(u_h)_h$
a subsequence, still denoted by $u_h$, such that on all subinterval
$[0,\tau']$ included in $\bigl[0,\tau\bigr]$, $u_h$
converges uniformly to $u$. In particular, thanks to
theorem~\ref{recap} we will have
\begin{equation*}
\begin{split}
&\lim_{h\to 0} \max
\{\vabs[U^m]{V^m}^2:0\le m\le \tau/h\} \\
&\quad\le \esssup \bigl\{ \vabs[u(t)]{\dot
u(t)}^2: t_0\le t \le t_0+\tau\bigr\},
\end{split}
\end{equation*} and for $h$ small enough we will have
\begin{equation*}
\begin{split}
&\lim_{h\to 0}\max
\{\vabs[U^m]{V^m}^2:0\le m\le \tau/h\} \\
&\quad\le C_2 =1+R
\max\{\vnorm{M(u)^{-1/2}}:\vabs{u-u_0}\le R\}.
\end{split}
\end{equation*}
Thanks to theorem~\ref{thr:9}, we can find $r_1$ such that for all
$\overline{u}\in K\cap B(u_0,R)$, for all $U^{l-1}$ and $U^l$ satisfying
condition $E(\overline{u}, r_1, C_2,h)$, it is possible to define a solution of
the scheme for $0\le (m-l)h \le \tau$, where $\tau$ is independent of
$h$. In particular, if we let $\tau'=\tau -\tau'/2$,
$l=\lfloor \tau'/h\rfloor$
and $\overline{u}=u(\tau')$, we can extend the scheme up to $m$
satisfying
\begin{equation*}
mh\le \min(lh +\tau'/2, \tau(R)),
\end{equation*}
which contradicts~\eqref{eq:47}. This proves the desired result.
\end{proof}

\bibliography{impact}
\bibliographystyle{siam}

\end{document}